\documentclass[11pt]{article}
\newcommand{\sm}[1]{\mbox{\small $#1$}}
\newcommand{\la}[1]{\mbox{\large $#1$}}
\newcommand{\La}[1]{\mbox{\Large $#1$}}
\newcommand{\LA}[1]{\mbox{\LARGE $#1$}}

\newcommand{\lnab}[1]{\la{\nabla}_{\!#1}}

\newcommand{\hlnab}[1]{\hat{\la{\nabla}}_{\!\!#1}}
\newcommand{\lnabp}[1]{\la{\nabla'}_{\!\!#1}}
\newcommand{\lnabo}[1]{\la{\nabla}^{\bot}_{\!\!#1}}
\newcommand{\lnabej}[1]{\la{\nabla}^{E_J}_{\!\!#1}} 
\newcommand{\lnabe}[1]{\la{\nabla}^{E}_{\!\!#1}}

\newcommand{\tlnab}[1]{\tilde{\la{\nabla}}_{\!\!#1}}

\newcommand{\qed}{\mbox{~~~\boldmath QED}}
\newcommand{\R}[1]{\mbox{\boldmath $R$}^{#1}}
\newcommand{\Co}{\mbox{\boldmath $C$}}
\newcommand{\Z}{\mbox{\boldmath $Z$}}
\newcommand{\ra}{\rightarrow}
\newcommand{\non}{\nonumber}
\newcommand{\ha}{\sm{\frac{1}{2}}}

\newcommand{\ddt}{\sm{\frac{d}{dt}}}
\newcommand{\al}{\alpha}
\newcommand{\be}{\beta}
\newcommand{\ga}{\gamma}
\newcommand{\bal}{\bar{\alpha}}

\newcommand{\Jw}{J_{\omega}}

\newcommand{\Ru}{\underline{R}}
\newtheorem{Lm}{Lemma}
\newtheorem{Pp}{Proposition}
\newtheorem{Th}{Theorem}
\newtheorem{Cr}{Corollary}

\newtheorem{Def}{Definition}

\begin{document}
\setcounter{page}{1}
\title{Transgression forms in dimension 4 \\[2mm] \small{by}\\[-5mm]}
\author{Isabel M.C.\ Salavessa$^{1*}$ and Ana Pereira do Vale$^{2**}$\\[2mm]}
\protect\footnotetext{$\!\!\!\!\!\!\!\!\!\!\!\!\!$
{\bf MSC 2000:} Primary: 53C42, 53C55, 53C25, 53C38; 
Secondary: 57R20, 57R45.\\
{\bf ~~Key Words:}
Characteristic classes, Transgression forms, Curvature tensors, 
almost complex structures, almost quaternionic K\"{a}hler,
Singular connections, dimension 4 and 8.\\
$*~\!\!$Partially supported by Funda\c{c}\~{a}o Ci\^{e}ncia e 
Tecnologia through POCTI/MAT/60671/2004 and Plurianual of CFIF.\\
$**$Partially supported by Funda\c{c}\~{a}o Ci\^{e}ncia e Tecnologia through
 POCTI/MAT/60671/2004.}

\date{\em Dedicated to Dmitri Alekseevsky on his 65th birthday\em ~~~~~~~~~
~~~~~~~~~~~~~~~~~~~~~~~~~~~~~~~~}
\maketitle ~~~\\[-9mm]
{\footnotesize $^1$ Centro de F\'{\i}sica das
Interac\c{c}\~{o}es Fundamentais, Instituto Superior T\'{e}cnico,
Edif\'{\i}cio Ci\^{e}ncia,\\[-1mm] Piso 3,
1049-001 LISBOA, Portugal;~~
e-mail: isabel@cartan.ist.utl.pt}\\
{\footnotesize $^2$ Centro de Matem\'{a}tica,
Universidade do Minho,
Campus de Gualtar\\[-1mm]
4710-057 BRAGA, Portugal;~~
e-mail: avale@math.uminho.pt}\\[4mm]
{\small {\bf Abstract:} 
We compute explicit transgression forms for the Euler and Pontrjagin 
classes of a Riemannian manifold  $M$ of dimension 4
 under a conformal change of the metric,
or a change to a Riemannian connection with torsion. 
These formulae describe
the singular set of some connections
with singularities on compact manifolds
as a residue formula  in terms of a polynomial of invariants.
We give some applications for minimal  submanifolds of
K\"{a}hler  manifolds.
We also
express the difference of the first Chern class of two  almost
complex structures, and in particular an obstruction 
to the existence of a homotopy between them, by a residue formula along 
the set of anti-complex points. 
Finally we take the first steps in the study  of obstructions
for two almost quaternionic-Hermitian  structures on a manifold of dimension 8
to have homotopic
fundamental forms or isomorphic twistor spaces. }
\markright{\sl\hfill  Salavessa--Pereira do Vale \hfill}
\section{Introduction}
\setcounter{Th}{0}
\setcounter{Pp}{0}
\setcounter{Cr} {0}
\setcounter{Lm} {0}
\setcounter{equation} {0}
\baselineskip .5cm
 Some  $k$-characteristic classes $Ch(E)$ of  a vector bundle $E$ over a 
manifold $M$, can be represented in the
 cohomology classes of $M$  in terms of the curvature tensor
defined with respect to a certain type of connections on $E$.
If $\lnab{}$, $\lnabp{}$ are 
such connections with curvature tensors
$R$ and $R'$, respectively, then
Chern-Weil theory states that 
\[Ch(R')= Ch(R)+ dT\]
where $T$ is a $(k-1)$-form on $M$. A way to specify such a $T$
is by pulling back each connection
 by the  projection $\pi:M\times [0,1]\ra M$,
 $\pi(p,t)=p$, and then take the connection 
$~\tlnab{}=t\lnabp{}+(1-t)\lnab{}~$ 
defined on $\pi^{-1}E$.
Denoting by $\tilde{R}$ its curvature tensor,
the Chern-Simons transgression $(k-1)$-form on $M$,  
 obtained by integration along [0,1]
 of the closed $k$-form $Ch(\tilde{R})$
on $M\times [0,1]$, 
\[T(\lnab{},\lnabp{}~)(X_1,\ldots, X_{k-1})=
\int_0^1 Ch(\tilde{R})(\ddt, X_1, \ldots, X_{k-1})dt
\]
satisfies $d\, T(\lnab{},\lnabp{}~) =Ch(R')- Ch(R)$.
To specify such $T$s is interesting by itself,
but also when $\lnabp{}$ is a connection with singularities.
This is the case when $(M,g)$ is a Riemannian $m$-manifold with
its Levi-Civita connection, and $(E, g_E, \lnabe{})$ is a Riemannian vector 
bundle of  rank $m$ such that there exists
a conformal bundle map $\Phi:TM \ra E$ which
vanishes along a singular set $\Sigma$. This bundle map induces on 
$M\!\sim\! \Sigma$ a  connection $\lnabp{}={\Phi^{-1}}^*\lnabe{}$, that makes 
$\Phi$ a parallel bundle map. 
This connection can be
seen as a singular Riemannian connection (with torsion) on $M$ 
with respect to a degenerated metric $\hat{g}$  on $M$, but $R'$
and 
$Ch(R')$ can be smoothly extended to $\Sigma$ by the identities
$ R'(X,Y,Z,W)= g_E(R^E(X,Y)\Phi(Z),\Phi(W))$, $Ch(R')\!=\!Ch(R^E)$. 
If $M$ is closed, and $Ch$ gives a  top rank form and an integral 
cohomology class, then 
$\int_M Ch(R)$ and $\int_M Ch(R')$ are finite  
 integers,  representing  invariants. Moreover, if the singular 
set $\Sigma$ is  sufficiently small and regular, the Stokes theorem  
reads $\int_{M\sim V_{\epsilon}(\Sigma)} dT=-\int_{\partial
V_{\epsilon}(\Sigma)}T$, 
where $V_{\epsilon}(\Sigma)$ is a tubular neighbourhood of $\Sigma$
of radius $\epsilon$,
and letting $\epsilon \ra 0$ may describe 
$Ch(E)\!-\!Ch(R)$ as a residue of $T$
 along $\Sigma$ and expressed in terms of the zeros of $\Phi$. 
This type of problem is studied 
in \cite{[H-L1]},\cite{[H-L2]},\cite{[H-Se]} , 
using currents. We  provide explicit formulae
of transgression forms for the cases of the Euler and Pontrjagin
classes. In
section 10  we give some  applications to minimal 4-submanifolds
in K\"{a}hler-Einstein manifolds.
In section 9 we describe an obstruction
for two almost complex structures on $M$ to be homotopic, measured
by the difference between their Chern classes, translated 
to  a residue formula on the set of anti-complex points. 
 In section 11  we introduce
 the study of obstructions for
two almost quaternionic-Hermitian structures on
a  Riemannian 8-manifold to have
isomorphic twistor spaces or homotopic fundamental forms.\\[3mm]
Let $(M,g)$ be an oriented  Riemannian manifold of dimension 4
with its Levi-Civita connection, curvature tensor $R^M$,
Ricci tensor $Ricci^M$ and scalar curvature $s^M$.
Our formulae are: 
\begin{Pp} Let  $f:M\ra \R{}$ be a smooth map, and $\hat{g}=e^f g$. Set
\begin{equation}
P(\nabla f)
=\La{(}2\Delta f+ \|\nabla f\|^2 - 2s^M\La{)}\nabla f
+4(Ricci^M)^{\sharp}(\nabla f)
-\nabla(\|\nabla f\|^2).
\end{equation}
Then\\[-10mm]
\begin{eqnarray}
{\cal X}(\hat{R}^M) &=& {\cal X}(R^M) +\frac{1}{32\pi^2}\,div_{g}
\LA{(} P(\nabla f)\LA{)}Vol_M, \\
p_1(\hat{R}^M)&=& p_1(R^M).\non \\[-3mm]\non
\end{eqnarray}
\end{Pp}
The invariance of the Pontrjagin class under a conformal
change of the metric is well known (see \cite{[At-B-P]} or remark 3). 
We thank Sergiu Moroianu for drawing our attention to this.
The authors do not recall to have seen formula (1.2) in the
literature: $div_{g}
\la{(} P(\nabla f)\la{)}$ is a 2nd-order differential operator on $\nabla f$.
If  $f\!=\!\log h$ where $h\!>\!0$ except at a finite set of 
 zeros $p_{\al}$ and   poles $p_{i}$ of  homogeneuos  order 
$2k_{\al}$ and $2k_i$ respectively ($k_{\al}, k_i>0$), 
then:\\[-3mm]
\[
\frac{1}{32\pi^2}\int_{M}div_g\La{(}P(\nabla\log h)\La{)}
Vol_{M}=\sm{\sum_{\al}}-\ha k^2_{\al}(k_{\al}+3)+ \sm{\sum_i}
-\ha k^2_i(-k_i+3).\\[3mm]\]
\begin{Th} Let $\Phi\!:\!TM\!\ra\! E$ be a conformal bundle map into a  
Riemannian  vector bundle $(E,g_E, \lnabe{})$ of rank 4 over $M$,
with coefficient of conformality given by a non-negative function
$h\!:\!M\!\ra\! \R{}$ with zero set $\Sigma$. Let $\hat{g}=hg$, and 
$S'\in C^{\infty}(\bigotimes ^2T^*M\otimes TM) $ defined away 
from ${\Sigma}$:\\[-3mm]
\[S'(X,Y)=\Phi^{-1}\La{(}\lnab{X}\Phi(Y)\La{)}-\ha d\log h(X)Y
-\ha d\log h(Y)X + \ha g(X,Y)\nabla \log h.\\[-2mm]\]
Then:\\[-7mm]
\begin{eqnarray*}
{\cal X}(R^E) &=& {\cal X}(R^M)-\frac{1}{4\pi^2}d\LA{(}
\langle {\cal S}'\wedge *\La{(}\hat{R}^M-\frac{1}{2}d{\cal S}'
-\frac{1}{3}({\cal S}')^2\La{)}\rangle_{\hat{g}}\LA{)}
+\frac{1}{32\pi^2}div_{g}
\LA{(} P(\nabla \log h)\LA{)}Vol_M\\
p_1(R^E) &=& p_1(R^M) -\frac{1}{2\pi^2}d\LA{(}
\langle {\cal S}'\wedge \La{(}\hat{R}^M-\frac{1}{2}d{\cal S}'
-\frac{1}{3}({\cal S}')^2\La{)}\rangle_{\hat{g}}\LA{)}\\[-5mm]
\end{eqnarray*}
where  $\hat{R}^M $ 
is the curvature tensor of $(M,\hat{g})$,
and ${\cal S}': TM\ra \bigwedge^2 TM$, 
$({\cal S}')^2: \bigwedge^2TM \ra \bigwedge^2TM$,  are defined by
\\[-6mm]
\begin{eqnarray*}
\langle{\cal S}'(X),Y\wedge Z\rangle_{\hat{g}}&=&\hat{g}(S'(X,Y),Z)
\\[-1mm]
\langle ({\cal S}')^2(X\wedge Y), Z\wedge W\rangle_{\hat{g}}&=&
\hat{g}(S'(X,Z), S'(Y,W))- \hat{g}(S'(X,W), S'(Y,Z))
\\[-1mm]
\langle \hat{R}^M(X\wedge Y),Z\wedge W\rangle_{\hat{g}}&=&
h\langle R^M(X\wedge Y) 
+ \phi\bullet g(X\wedge Y), Z\wedge W\rangle_g\\[-9mm]
\end{eqnarray*}
\[
\phi=\frac{1}{2}\La{(}\sm{-\frac{\|\nabla \log h\|^{2}}{4}}g
+\ha d\log h \otimes d\log h- Hess(\log h)\La{)}.\]
Moreover, if $\, d\Phi=0$, then ${\cal X}(R^E)={\cal X}(R^M)+\frac{1}{32\pi^2}
div_g(P(\nabla \log h))Vol_M$ and
 $p_1(R^E)=p_1(R^M)$.\\[-3mm]
\end{Th}

The angle $\theta\in [0,\pi]$ between two positive $g$-orthogonal 
almost complex
structures $J_0$  and $J_1$ on $M$ is defined by
$ \cos\theta =\frac{1}{4}\langle J_0,J_1\rangle $, and
 $\Sigma=\{p\in M: J_1(p)\!=\!-J_0(p)\}=\cos\theta^{-1}(-1)$ is  the
set of \em anti-complex points. \em If $J_1$ is generic, $\Sigma$ is
a surface of $M$ and an orientation can be given. 
Under the usual identification of $J_0$ with its
K\"{a}hler form $\omega_0$,
 the orthogonal complement $E_{J_{0}}$ of $\R{}J_{0}$ in
$\bigwedge^2_+TM$, is a complex line bundle over $M$
with complex structure  $"J_0"$.
We denote by $\tilde{H}(p)$ the orthogonal projection of $J_1(p)$
into $E_{J_0}$.
Let $N^1\Sigma$ be the total set of the unit normal bundle
of $\Sigma$ and $d_{N^1\Sigma}$ its Lebesgue measure.
For each $(p,u)\in N^1\Sigma$ define $\kappa(p,u)\geq 1$ 
the order of the zero of $\phi_{(p,u)}(r)= (1+\cos\theta)(exp_{p}(ru))$ 
at $r=0$.  We say that $(1+\cos\theta)$
has a \em controlled zero set \em if  
there exist a non-negative integrable function $f\!:\!N^1\Sigma
\ra [0,+\infty]$ and $r_0>0$ s.t. 
$~sup_{0<r<r_0}|r\frac{d}{d r}\log(\phi_{(p,u)}(r))|\leq f(p,u)~$
a.e.\ $(p,u)\in N^1\Sigma$. For example, this holds if for  all $(p,u),$
 $\phi_{(p,u)}(r)$ is a polynomial  function on $r$ with 
coefficient of lowest order uniformly bounded away from zero.
 Some weaker conditions can be given on
$\frac{d^k}{d r^k}\phi_{(p,u)}(0)$, $\frac{d^{k+1}}{d r^{k+1}}
\phi_{(p,u)}(r)$ and $\frac{d^{k+2}}{d r^{k+2}}
\phi_{(p,u)}(r)$, where $k=\kappa(p,u)$, 
to guarantee controlled zero set (see Prop.\ 9.4).
For each $p\in \Sigma$, $S(p,1)$ denotes the
unit sphere of $T_p\Sigma^{\bot}$ and $\sigma_{d'}$ its volume. The function
 $\tilde{\kappa}(p)=\frac{1}{\sigma_{d'}}
\int_{S(p,1)}\kappa(p,u)d_{S(p,1)} u$  is the \em average order \em
of the zero $p$ of $(1+\cos\theta)$, in the normal direction.
Assume $M$ is  compact.
\\[-4mm]
\begin{Th} If  $J_0$ is almost K\"{a}hler  
then for any almost complex structure $J_1=\cos\theta J_0+\tilde{H}$
\\[-2mm]
\begin{equation}
\int_M\La{(}c_1(M,J_1)-c_1(M,J_0)\La{)}\wedge \omega_0\
= \frac{1}{4\pi}\int_M div((\tilde{T}J_0)^{\sharp})Vol_M\\[1mm]
\end{equation}
where  $\tilde{T}$ is the 1-form on $M\!\sim\! \Sigma$,
$\tilde{T}(X)=\frac{1}{(1+\cos\theta)}\langle\lnab{X}\tilde{H},
J_0\tilde{H}\rangle.$
In the particular case that  $\lnab{}^{E_{J_0}}\tilde{H}$ is 
$J_0$-anti-complex we have\\[-2mm]
\begin{equation}
(1.3)= -\frac{1}{4\pi}\int_M \Delta \log (1+\cos\theta)Vol_M.
\end{equation}
In this case, assume $\Sigma$ is a finite disjoint 
union of closed oriented submanifolds $\Sigma^i$ of dimension 
$d_i\!\leq\! 2$. Let 
$\bigcup_{\ga}k_{i\ga}$ be the range set of $\kappa$ on $N^1\Sigma^i$
 and  let $N^1\Sigma^i_{\ga}=\kappa^{-1}(k_{i\ga})$. 
If $\kappa$ is bounded a.e.
and $(1+\cos\theta)$ has a controlled zero set
then a residue formula along $\Sigma$ is obtained:\\[-2mm]
\begin{equation}
(1.3)
=\ha\sum_{i: d_i=2}\int_{\Sigma^i}\tilde{\kappa}(p)Vol_{\Sigma^i}= 
\ha\sum_{i: d_i=2}\sum_{\ga}k_{i\ga} \, d_{N^1\Sigma^i}(N^1\Sigma_{\ga}^i).
\end{equation}
and so\\[-10mm]
\[~~~
\int_M\langle c_1(M,J_1), \omega_0\rangle Vol_M \geq 
\int_M\langle c_1(M,J_0), \omega_0\rangle Vol_M =\frac{1}{4\pi}
\int_M (s^M +\ha \|\lnab{}\omega_0\|^2 )Vol_M,\\[-1mm]
\]
with equality iff $d_i\leq 1$ $\forall i$. Thus,  if 
$c_1(M,J_1)=c_1(M,J_0)$ (in $H^2(M,\R{})$)  
$\Sigma_i=\emptyset$ for $d_i= 2$.
\end{Th}
\em  Remark $1$. \em  Pairs of almost complex structures with or without the
same properties may exist on a  manifold.
 Alekseevsky [1] discovered examples of 
simply-connected non-compact Riemannian manifolds admitting 
a non-integrable almost-K\"{a}hler structure and also an integrable  
non-K\"{a}hler complex  structure,  namely some solvable groups of dimension 
$4(4+p+q)$ $ (q\neq 0)$.\\[2mm]
We also prove in section 9:
\begin{Pp} If $M$ is  compact, $J_0$ is  K\"{a}hler
and $J_1$ is almost K\"{a}hler
then $\theta$ is constant. Thus, if $\cos\theta\neq \pm 1$,  $J_0$ and $J_1$
are homotopic and  define a hyper-K\"{a}hler structure on $M$.
\end{Pp}
If $(M^8,g)$ is an oriented 8-dimensional manifold and $Q_0$ and $Q_1$
are two almost quaternionic Hermitian structures (see [2]
for definitions) we define an angle
$\theta\in[0,\pi]$, by $\cos\theta =\frac{3}{10}\langle
\Omega_0,\Omega_1\rangle, $ where $\Omega_i$ are the corresponding
fundamental 4-forms. Let $E_i$ be the corresponding rank 3 vector bundle
generated by the twistor space of $Q_i$. In section 11 we prove:
\begin{Pp} If $Q_0$ is quaternionic K\"{a}hler and  $Q_1$ 
is almost quaternionic K\"{a}hler and $M$ is compact, then
$\theta$ is constant.  If $\cos\theta\neq -1$ then $\Omega_0$ and 
$\Omega_1$ are homotopic 4-forms in $H^4_+(M;\R{})$, and  if
$Q_1$ is also quaternionic K\"{a}hler,  $p_1(E_0)=p_1(E_1)$.
Furthermore, in the later case, $(a), (b)$ or $(c)$ must hold: 
$(a)$ $E_0=E_1$;
$(b)$ $E_0\cap E_1$ has rank one, $M$ is K\"{a}hler and both $Q_0, Q_1$ 
are locally hyper-K\"{a}hler structures; 
$(c)$ $E_0\cap E_1=\{0\}$.
\end{Pp}
We observe that the problem on a compact 4-manifold "almost K\"{a}hler 
+ Einstein implies K\"{a}hler", also called the Goldberg-conjecture,
 is not completely solved. Salamon in \cite{[Sa]} gives an example of a
compact 8-dimensional almost quaternionic K\"{a}hler manifold that
is not quaternionic K\"{a}hler. We thank the referee for drawing our
attention to this reference.
\section{Curvature tensors in dimension 4}
\setcounter{Th}{0}
\setcounter{Pp}{0}
\setcounter{Cr} {0}
\setcounter{Lm} {0}
\setcounter{equation} {0}
\baselineskip .5cm
Let $V$ be a vector space of dimension $4$ and with a inner product $g$. 
We identify $\bigwedge^2:=\bigwedge^2V$ with  
$Skew(V)$ (the space of  skew-symmetric
endomorphisms of $V$) and with $\bigwedge^2V^*$ in the standard way,
considering $Skew(V)$ as a subset of $V^*\otimes V$ with half
of its usual Hilbert-Schmidt inner product. 
Consider the vector subspaces
${\cal R}$ and ${\cal R}^{\bot}$ of
the symmetric and  skew-symmetric linear endomorphisms of 
$\bigwedge^2$,  respectively. 
${\cal R}$, or more generally $L(\bigwedge^2;L(V,V))$, is defined
as the space of curvature tensors of $V$ (see \cite{[Be1]}, \cite{[Be2]},
or  \cite{[Si-T]} for details). We recall some definitions.

If $R \in L(\bigwedge^2;L(V,V))$
let $\underline{R}\in L(V;L(\bigwedge^2;V))$ given by
$\underline{R}(Z)(X\wedge Y):=R(X\wedge Y)Z$. 
We also use the following notation:
$R(X,Y,Z,W)=g(R(X\wedge Y)Z,W)=
\underline{R}(Z,W,X,Y)$.
If $R \in L(\bigwedge^2;\bigwedge^2)$, $\Ru=R^T$ (transposed).
We  assume that $V$ is with a given orientation, and so
the star operator $*\in L(\bigwedge^2;\bigwedge^2)$
 splits $\bigwedge^2$ into its eigenspaces $\bigwedge^2_{\pm}$,
 corresponding to the eigenvalues $\pm 1$, defining respectively
the space of selfdual and of anti-self-dual two forms. Then 
$R\in L(\bigwedge^2;\bigwedge^2)$ splits as 
$ R= R^+_+\oplus R^+_-\oplus R^-_+\oplus R^-_-$ with
$R^+_{\pm}\in L(\bigwedge^2_{\pm};\bigwedge^2_{+})$,
$R^-_{\pm}\in L(\bigwedge^2_{\pm};\bigwedge^2_{-})$.
The Ricci tensor  $Ricci:L(\bigwedge^2;\bigwedge^2)\ra L(V,V)$
and the scalar curvature are given respectively  by
$g( Ricci(R)(X),Y)= Ricci_R(X,Y)= tr\{Z\ra R(X,Z)(Y)\}$
and $ s_R=tr\,(Ricci_R)=2\langle R, Id\rangle = 2 tr\,(R)$.
The sectional curvature  of $R\in L(\bigwedge^2;\bigwedge^2)$ 
is  denoted by $\sigma_R(P)=R(X,Y,X,Y)$ for each 2-plane $P$ spanned 
by an o.n.b.\ $\{X,Y\}$, and 
the  Bianchi map $b:L(\bigwedge^2; L(V;V))\ra \bigwedge^3V^*\otimes V
\subset L(\bigwedge^2; L(V;V))$
is defined by 
$~b(R)(X_1,X_2)(X_3)= R(X_1,X_2)(X_3)+ R(X_3,X_1)(X_2)+ R(X_2,X_3)(X_1).$
Then 
$b:{\cal R}\ra {\cal R}$,  ${\cal R}^{\bot}\ra {\cal R}^{\bot}$.
If $R\in {\cal R}$, $b(R)=\ha tr(*R)*.$
Let $t(R)=\frac{1}{6}\langle b(R),*\rangle$. We denote by $Sym(V)$ the
space of symmetric  endomorphisms of $V$ and by $Sym_0(V)$
its subspace of trace-free endomorphisms.
\\[3mm]
A complex 2-plane $\tau$ of $V^c$ is said to be
totally isotropic (t.i.) if 
$g(v,v)=0$ $\forall v\in \tau$. 
If $(e_i)$ is  an o.n. basis  and
$\al=e_1\!+\! ie_2$, $\be=e_3\!+\!ie_4$, then $\tau=span_{\Co}\{\al, \be\}$
defines a t.i.\  complex plane. Set $R(ijkl)=R(e_i,e_j,e_k,e_l)$. 
The isotropic sectional curvature w.r.t $R\in {\cal R}$  at $\tau$
is given by (see \cite{[Mi-Mo]})
\begin{equation}
K_{isot}(R)(\tau)=\frac{R(\al\wedge \be,\bar{\al}\wedge \bar{\be})}{
\|\al\wedge \be\|^2}= \frac{1}{4}\LA{(}
R(1313)+R(1414)+R(2323)+R(2424)-2R(1234)\LA{)}.
\end{equation}
The \em Kulkarni-Nomizu product \em of
 $\xi,\phi\in \bigotimes^2V^*$  is a symmetric product defined by
\[\xi\bullet \phi (X,Y,Z,W)=\xi(X,Z)\phi(Y,W)+\xi(Y,W)\phi(X,Z)
-\xi(Y,Z)\phi(X,W) -\xi(X,W)\phi(Y,Z).\]
We have $Id=Id_{\bigwedge^2}=\ha g\bullet g$.
If $R\in L(\bigwedge^2;L(V;V))$ it is  defined  the 
\em Weitzenb\"{o}ck  operator \em  $A(R)$
\[ \langle A(R)(X,Y),Z\wedge W\rangle:= 
Ricci_R\bullet g(X,Y,Z,W)+ 2R(Z,X,Y,W)-2R(W,X,Y,Z).\]
If $R\in {\cal R}$ then $A(R)\in {\cal R}$
and $A(R)=Ricci_R\bullet g -2R + 2b(R)$. 
In this case, $Ricci_{A(R)}=s_R g$, $b(A(R))=4b(R)$.
Furthermore, if $b(R)=0$,
$A(R):\bigwedge^{\pm}\ra \bigwedge^{\pm}$. Note that $*\in {\cal R}$,
and if $R\in {\cal R}$, then $*R+R*,*R*$ still lie in ${\cal R}$ and
$*R-R*\in {\cal R}^{\bot}$. 
Straightforward computations shows:
\begin{Pp} Let $R\in {\cal R}$ and $\tau=span_{\Co}\{e_1+ie_2,
e_3+ie_4\}$  with $e_1,\ldots,e_4$ an
o.n. basis with orientation $\epsilon$. $K_{isot}$ is 
 computed at $\tau$. Then
\[\begin{array}{llll}
\sm{Ricci_{Id}=3 g} &\sm{s_{Id}= 12} & \sm{b(Id)= 0}& 
\sm{K_{isot}(Id)=1}\\  
\sm{Ricci_{*\mbox{}}= 0} &\sm{ s_{*} = 0}& \sm{b(*)= 3*}&
\sm{K_{isot}(*)=-\ha \epsilon}\\
\sm{Ricci_{*R*}= \ha s_R\, g-Ricci_R}& \sm{s_{*R*}= s_R}
& \sm{b(*R*)=b(R)}& \sm{K_{isot}(*R*)=K_{isot}(R)}\\
\sm{Ricci_{R*}= Ricci_{*R}=t(R)g} & \sm{s_{*R}=s_{R*}= 4t(R)}
& \sm{b(*R\!\!-\!\!R*)= 0} & \sm{K_{isot}(*R\!\!-\!\!R*)
=2(\sigma_R(34)\!\!-\!\!\sigma_R(12)).}
\end{array}\]
\end{Pp}
For $R, Q \in L(\bigwedge^2;\bigwedge^2)$, and $X,Y,Z,W\in V$
we use the following notation:
$R(X\wedge Y)\wedge Q(Z\wedge W)=$ $
\langle R(X\wedge Y), (*Q)(Z\wedge W)\rangle Vol_V
~\in \bigwedge^4 V^*$ 
\begin{Def} 
 The  Euler form ${\cal X}(R)$ and the Pontrjagin form $p_1(R)$ 
of $R\in L(\bigwedge^2;\bigwedge^2)$ are the 4-forms:\\[-7mm]
\begin{eqnarray}
4\pi^2{\cal X}(R) &=& { \Ru(e_1\wedge e_2)\wedge \Ru(e_3 \wedge e_4)-
\Ru(e_1\wedge e_3)\wedge \Ru(e_2 \wedge e_4)+
\Ru(e_1\wedge e_4)\wedge R(e_2\wedge e_3)}\non\\
&=&\ha \langle R, * R *\rangle Vol_V
=\ha \langle *R,R*\rangle Vol_V\\[2mm]
4\pi^2{p}_1(R)&=&\Ru(e_1\wedge e_2)\wedge \Ru(e_1\wedge e_2)
+\Ru(e_1\wedge e_3)\wedge \Ru(e_1\wedge e_3)
+\Ru(e_1\wedge e_4)\wedge \Ru(e_1\wedge e_4)\non\\[-1mm]
&& + \Ru(e_2\wedge e_3)\wedge \Ru(e_2\wedge e_3) +
\Ru(e_2\wedge e_4)\wedge \Ru(e_2\wedge e_4)
+ \Ru(e_3\wedge e_4)\wedge \Ru(e_3\wedge e_4)\non\\
&=& \langle R^T, * R^T\rangle Vol_V =\langle R, R*\rangle Vol_V
\end{eqnarray}
\end{Def}
Note that
$\langle R,Q\rangle\! =\! \langle R^T, Q^T\rangle$. Thus  
 ${\cal X}(R^T)={\cal X}(R)$,  and $
p_1(R^T)\!-\!p_1(R)=\langle R, *R\!-\!R*\rangle.$\\[4mm]
A positive $g$-orthogonal complex structure $J$ on $V$ is
a complex structure that induces the orientation of $V$ and
it is a linear isometry. Such structures are in 1-1
correspondence with the  elements $\omega_J$ of $\bigwedge^2_+V$ of
norm $\sqrt{2}$, by $\omega_J(X,Y)= g(JX,Y)$. The condition of orthogonality
between two of such $J$ is equivalent to the anti-commuting condition,
for the multiplication of complex structures corresponds to the quaternionic
multiplication of two unit pure imaginary vectors of $\R{4}$.
For each $J$ and $R\in L(\bigwedge^2;\bigwedge^2)$  
we  define the 2-forms on $V$
\begin{eqnarray}
Ricci_{J,R}(X,Y)&:=&\langle R(X\wedge Y),\omega_J\rangle=
\underline{R}(\omega_J)(X,Y)~~~~~~~\\
\Psi_J(R)(X,Y) &:=& -\ha Ricci_R(X,JY)+\ha Ricci_R(Y,JX)
\end{eqnarray}
If $R\in {\cal R}$, then $\Psi_J(R)(X,Y)\!=\! Ricci^{(1,1)}(JX,Y)$.
Usually $Ricci_{J,R}$ is denoted by $Ricci_*$ and named by
\em star-Ricci \em  but we do
not use that notation to avoid confusion with the Ricci of 
the  curvature tensor $*$. The \em star-scalar \em curvature is
 $s_{J,R}\!=2\langle Ricci_{J,R},\omega_J\rangle$, and 
${b}_J(R)\in \bigwedge^2V^*$ is defined by
$~{b}_J(R)(X, Y) = -\ha tr(Z\ra Jb(R)(X,Y)(Z)).$ 
If $P\in L_J(V;V)$, 
that is $P\circ J= J\circ P$ then both the complex trace $tr_J(P)$
and the  complex determinant $ det_{J}(P)$ are well defined. 
If $X_1,Y_1=JX_1,X_2,Y_2=JX_2$ is a real basis of $V$
and  for $\al= 1,2$ denote $"\al" := W_{\al}=\ha(X_{\al}-iY_{\al})$, 
$\bal:= W_{\bar{\al}}=\overline{W_{\al}}$, then
$\sm{tr_{J}(P)}=\sm{\sum_{\al}W^{\al}
_* P^c(W_{\al})}$,
$\sm{det_J(P)}=\sm{   det\, 
[W^{\al}_*P^c(W_{\be})]}$
where $(W_*^{\al}, W_*^{\overline{\al}}=\overline{W_*^{\al}})$
is the complex dual basis of  $\{W_{\al},W_{\bar{\al}}\}$,
 and $P^c$ denotes the complex 
linear extension of  $P$ to $V^c$. 
If $P\in L_J(V;V)\cap \bigwedge^2$, then
$ \langle P,\omega_J\rangle  =tr_J(P).$
If  $R\in L(\bigwedge^2;L_J(V;V))$, then $R$ is said to be
\em $J$-invariant. \em  \\[-3mm]
\begin{Def} If $(V,J,g)$ is Hermitian, the \em first \em and \em 
second Chern form \em  of $R\in L(\bigwedge^2;L_J(V;V))$  w.r.t. $J$,
are respectively\\[-3mm]
\begin{equation}
 c_1(R,J) = -\frac{i}{2\pi} Tr_J(\underline{R})~~~~~~~~~~
 c_2(R,J) = -\frac{1}{4\pi^2}det_J(\underline{R}).\\[1mm]
\end{equation}
\end{Def}
It follows that, if $R\in L(\bigwedge^2;L_J(V;V)\cap \bigwedge^2)$
then $c_1(R,J)={\frac{1}{2\pi}} Ricci_{J,R}$, 
$~c_2(R,J)=  {\cal X}(R)$, 
$ p_1(R) = c_1(R,J)\wedge c_1(R,J)-2c_2(R,J)$, and
$b_J(R) = Ricci_{J,R}- \Psi_J$.
If $R\in{\cal R}$ and is $J$-invariant, then 
$Ricci_R(JX,JY)=Ricci_R(X,Y)$, and  
$\Psi_J(R)(X,Y)= Ricci_R(JX,Y)$.\\[5mm] 
We denote by $E_J$ the rank 2 subspace of
$\bigwedge^2_+V$   defined by the orthogonal complement of
$\R{}\{\omega_J\}$\\[-3mm]
\[\sm{\bigwedge^2_+}V=\R{}\{\omega_J\}\oplus E_J, ~~~~~\mbox{where}~~~
 E_J=\{\omega\in \sm{\bigwedge^2_+}V: \omega(JX,JY)=-\omega(X,Y)\}\]
and  a canonic complex structure can be given to $E_J$:
$\tilde{J}\omega(X,Y)=-\omega(JX,Y)$. Let $\{e_1,e_2,e_3,e_4\}$ be a
d.o.n.b.\ of $V$, giving a corresponding o.n.b.\ 
$\sqrt{2}\omega_{\sigma}$ of $\bigwedge^2_+V$,
 $\omega_1=e_1\wedge e_2 +e_3\wedge e_4$,
$\omega_2=e_1\wedge e_3-e_2\wedge e_4$, 
$\omega_3=e_1\wedge e_4+e_2\wedge e_3$, and
let  $J_{\sigma}$ be defined by $\omega_{\sigma}=g(J_{\sigma}(\cdot), \cdot)$. 
The 2-forms  $\omega_{2},\omega_{3}$ span 
$E_J$, and $\tilde{J}\omega_{2}=\omega_{3}$ corresponds to $J_1J_2=J_3$.
Any such o.n. (of norm $\sqrt{2}$) basis $(\omega_1,\omega_2,\omega_3)$,
where $\omega_3="\omega_1\omega_2"$ defines a canonic orientation
on $\bigwedge^2_+V$. 
\section{Almost Hermitian 4-manifolds}
\setcounter{Th}{0}
\setcounter{Pp}{0}
\setcounter{Cr} {0}
\setcounter{Lm} {0}
\setcounter{equation} {0}
\baselineskip .5cm
Assume $(M,J,g)$ is an almost Hermitian 4-manifold with its Levi-Civita
connection $\lnab{}$, and we use the above notation taking $V=T_pM$.
 $(M,J,g)$ is
said to be \em almost K\"{a}hler \em 
if the K\"{a}hler form $\omega_J(X,Y)=g(JX,Y)$
is closed, and K\"{a}hler if $\omega_J$ is parallel.
The latter  is equivalent to $J$ to be  almost K\"{a}hler
and integrable.  Since $J^2=-Id$ and $J$ is $g$-orthogonal,  
$\forall X,Y,Z\in T_pM$
\begin{equation}
\lnab{Z}J(JX)=-J\lnab{Z}J(X), ~~~~~~
g(\lnab{Z}J(X),Y)=-g(X, \lnab{Z}J(Y)).
\end{equation}
We use the following sign for curvature tensors
$R(X,Y)=-\lnab{X}\lnab{Y}+\lnab{Y}\lnab{X} +\lnab{[X,Y]}$,
$~\forall X,Y\in T_pM,$ and
 denote by $R^M$ the curvature tensor of $M$, 
$Ricci^M_J=Ricci_{J,R^M}$~, 
$s^M_J=s_{J,R^M}$, and by $\langle\cdot,\cdot\rangle$ the Hilbert
Schmidt inner product on $\bigwedge^kTM^*$.
The   Weitzenb\"{o}ck formulae for $\omega_J$ read (see e.g. \cite{[Be2]})
\begin{eqnarray}
\Delta \omega_J&=&-tr \lnab{}^2\omega_J+ A(R^M)(w_J) \\
0 = \ha\Delta^+\|\omega_J\|^2 &=&-\langle \Delta \omega_J,\omega_J\rangle
+\|\lnab{}\omega_J\|^2 +\langle A(R^M)\omega_J, \omega_J\rangle.
\end{eqnarray}
Since $\omega_J$ is a self-dual form, 
$\|d\omega_J\|=\|\delta \omega_J\|$. Thus $J$ 
is almost K\"{a}hler iff $\delta \omega_J=0$. In this case
$\omega_J$ is harmonic.
If $J$ is not K\"{a}hler we may use the
 \em canonical Hermitian  \em connection (see \cite{[Gau]})
\[\tlnab{X}Y =\lnab{X}Y-\ha J(\lnab{X}J (Y)).\]
This connection satisfies $\tlnab{}g=\tlnab{}J=0$, but has
torsion $\tilde{T}(X,Y)=-\ha JdJ(X,Y)$.
This is a $U(2)$-connection on $M$, and so its curvature tensor
$\tilde{R}$ is $J$-invariant.
The Chern classes of $M$ using $\tilde{R}$ satisfy
$2\pi c_1(M,J)= {R}icci_{J,\tilde{R}} $, 
and a direct computation
shows that
\[{R}icci_{J,\tilde{R}}(X,Y)= Ricci_J^M(X,Y)+\eta_J(X,Y),~~~~~~~~~~~ 
s_{J,\tilde{R}}=s^M_J +2\langle \eta_J, \omega_J\rangle\]
where $\eta_J$ is the 2-form on $M$
\begin{equation}
\eta_J(X,Y) = \sm{\frac{1}{4}}\langle J\lnab{X}J, \lnab{Y}J\rangle, 
\end{equation}
and the inner product is the Hilbert-Schmidt inner product on
$TM^*\otimes TM$.
If $(M,J,g)$ is almost K\"{a}hler then (see \cite{[Ap-D]} for a survey)
\begin{equation}
\lnab{JX}J= -J(\lnab{X}J).
\end{equation}
In this case
$\tilde{T}^{(1,1)}=0$, and
$\langle \eta_J,\omega_J\rangle=
-\frac{1}{8}\|\lnab{}J\|^2.$ 
Furthermore, for almost K\"{a}hler $J$ (see \cite{[Bl]})
\begin{eqnarray}
s^M_J-s^M &=&\|\lnab{}\omega_J\|^2\\
4\pi \langle c_1(M,J), \omega_J\rangle &=&
2\langle \tilde{R}(\omega_J), \omega_J\rangle = s^M
+\ha\|\lnab{}\omega_J\|^2=\ha(s^M_J+s^M).
\end{eqnarray}
We consider on $E_J$ the induced connection $\lnabej{}$ from the connection
$\lnab{}^{+}$ of $\bigwedge^2_+$.
\begin{Pp} If $(M,J,g)$ is almost Hermitian 
\begin{eqnarray*}
2\pi c_1(E_J)(X,Y)&=&
Ricci_{J}^M(X,Y)+\eta_J(X,Y)=2\pi c_1(M,J)\\[-1mm]
p_1(\sm{\bigwedge^2_+}T_pM)&=&c_1(E_J)\wedge c_1(E_J).
\end{eqnarray*}
\end{Pp}
\em Proof. \em   Let $\frac{\omega_2}{\sqrt{2}}$, 
$\frac{\omega_3}{\sqrt{2}}$ be a local d.o.n. frame   of $E_J$,
and  $(\cdot )^{E_J}$ denote
 the orthogonal projection onto $E_J$. 
Since $\|\omega_J\|=\sqrt{2}$ is constant,
$ \lnab{Y}^+\omega_J$ is a section of $E_J$. Now,
\begin{eqnarray*}
\lnabej{Y}\lnabej{X}s &=&  (\lnab{Y}^+\lnabej{X}s)^{E_J}
= (\lnab{Y}^+\lnab{X}^+s-\ha\langle \lnab{X}^+s, \omega_J\rangle
\lnab{Y}^+\omega_J)^{E_J}\\
&=& (\lnab{Y}^+\lnab{X}^+s)^{E_J}+\ha\langle s, \lnab{X}^+\omega_J\rangle
\lnab{Y}^+\omega_J.
\end{eqnarray*}
Thus
\begin{equation}
R^{E_J}(X,Y)s = (R^+(X,Y)s)^{E}+\ha\langle s, \lnab{X}^+\omega_J\rangle
\lnab{Y}^+\omega_J - \ha\langle s, \lnab{Y}^+\omega_J\rangle
\lnab{X}^+\omega_J.
\end{equation}
The curvature tensor of  $\bigwedge^2_+$  satisfies
$\langle R^+(X,Y)\omega_2, \omega_3\rangle = 2R^M(X,Y)\omega_J.$
$E_J$ is a complex line bundle over $M$ and so it has a (real) volume
element $Vol_{E_J}$.
Using (3.1) we easily see that
another way to express $\eta_J$ is
$\eta_J(X,Y)=\ha Vol_{E_J}(\lnab{X}^+\omega_J,
\lnab{Y}^+\omega_J).$
Then
\begin{equation}
2\pi c_1(E_J)(X,Y):=\langle R^{E_J}(X,Y)\frac{\omega_2}{\sqrt{2}},
\frac{\omega_3}{\sqrt{2}}\rangle=  R^M(X,Y)\omega_J+ \eta_J(X,Y).
\end{equation}
Since  $c_2(\Co\{\omega_J\}\oplus E^{c}_J)
=c_1(\Co\{\omega_J\})\wedge c_1(E^{c}_J)-c_2(E^c_J)$ 
and $c_1(E^{c}_J)=0$
we have $p_1(\bigwedge^2_+TM)=-c_2((\bigwedge^2_+TM)^c)
=-c_2(\Co\{\omega_J\}\oplus E^{c}_J)
= c_1(E_J)\wedge c_1(E_J).$ \qed \\[2mm]
\em Remark 2. \em 
The \em  canonical line bundle \em  w.r.t $J$ is  
the bundle ${\cal K}_J=\bigwedge^2(TM^{(1,0)})^*
=\bigwedge^{(2,0)} $
of the complex 2-forms of type $(2,0)$, and one has 
$(\bigwedge^2_+TM)^c =\Co(\omega_J)\oplus \bigwedge^{(2,0)}
\oplus \bigwedge^{(0,2)}$.
The bundle $E_J$ is isomorphic to  the realification of the
anti-canonical bundle ${\cal K}^{-1}_J=\bigwedge^{(0,2)}$. 
As a complex line bundle,
$\omega\ra  (\omega_c)^{(0,2)}$ is the complex isomorphism.
Therefore, $c_1(E_J)=c_1({\cal K}^{-1}_J)=-c_1({\cal K}_J)$.
An almost complex structure defines 
a canonic spin-c structure ${\bf s}=P_{Spin^c}M$ on $M$ with 
canonic line bundle $\bigwedge^{(2,0)}$ and 
complex spinor bundle  $S\!=\!\bigwedge^{(0,*)}\!=\! S^+\oplus S^-$, 
$S^+\!=\!\bigwedge^{(0,0)}\!+\!\bigwedge^{(0,2)}$,
$S^-\!=\!\bigwedge^{(0,1)}$, where $TM$ acts by Clifford multiplication.
The Chern class of ${\bf s}$ is also given by $c_1({\cal K}_J)$.\\[3mm]
Finally we observe the following:
\begin{Lm} If $J_1$ and $J_2$ are two anti-commuting $g$-orthogonal
almost K\"{a}hler complex structures on $M$, then $J_3=J_1J_2$ is also
 almost K\"{a}hler iff $g(\lnab{Z}J_2(X),J_1X)=0$ $\forall X\in T_pM$.
That is the case if either $J_1$ or $J_2$ is K\"{a}hler.
In that case  $(J_1,J_2,J_3)$ is in fact   an hyper-K\"{a}hler
structure.
\end{Lm}
\em Proof. \em  It is sufficient to find an equivalent
condition for   $\delta \omega_3\!=\!0$. Using (3.1), (3.5),
 $\forall X\in T_pM$, 
\begin{eqnarray*}
\lefteqn{\lnab{J_3X}J_3(J_3X)+\lnab{X}J_3(X)=
\lnab{J_1J_2X}(J_1J_2)(J_3X)+\lnab{X}(J_1J_2)(X)=}\\
&=&\lnab{J_1J_2X}J_1(J_2J_3X)+ J_1(\lnab{J_1J_2X}J_2(J_3X))
+\lnab{X}J_1(J_2X) +J_1(\lnab{X}J_2(X))\\
&=&-J_1(\lnab{J_2X}J_1(J_2J_3X))-J_1(\lnab{J_2J_1X}J_2(J_3X))
+\lnab{X}J_1(J_2X) +J_1(\lnab{X}J_2(X))\\
&=&\lnab{J_2X}J_1(J_1J_2J_3X)+J_1J_2(\lnab{J_1X}J_2(J_3X))
+\lnab{X}J_1(J_2X) +J_1(\lnab{X}J_2(X))\\
&=&-\lnab{J_2X}J_1 (X)+\lnab{X}J_1(J_2X)
-J_1\lnab{J_1X}J_2(J_2J_3X)+J_1(\lnab{X}J_2(X))\\
&=&\LA{(}-\lnab{J_2X}J_1 (X)+\lnab{X}J_1(J_2X)\LA{)}
+J_1\LA{(}-\lnab{J_1X}J_2(J_1X)+\lnab{X}J_2(X)\LA{)}
\end{eqnarray*}
Using (3.5) from $d\omega_1(X,J_2X,Z)=0$ $ \forall X,Z$, we have
$g(-\lnab{J_2X}J_1 (X)+\lnab{X}J_1(J_2X),Z)=-g(\lnab{Z}J_1(X),J_2X)$.
If $X$ is a unit vector of $T_pM$ then $X,J_1X,J_2X,J_3X$ is an o.n.b.
Hence
\begin{eqnarray*}
\lefteqn{-g(\delta J_3,Z)= g\La{(}\lnab{X}J_3(X)+\lnab{J_3X}J_3(J_3X)
+\lnab{J_1X}J_3(J_1X)+\lnab{J_3J_1X}J_3(J_3J_1X), Z\La{)}=}\\
&=&-g(\lnab{Z}J_1(X),J_2X)-g(-\lnab{J_1X}J_2(J_1X)+\lnab{X}J_2(X), J_1Z)\\
&&-g(\lnab{Z}J_1(J_1X),J_2J_1X)
-g(-\lnab{X}J_2(X)+\lnab{J_1X}J_2(J_1X), J_1Z)\\
&=&-g(\lnab{Z}J_1(X),J_2X)-g(J_1\lnab{Z}J_1(X),J_1J_2X)
=-2g(\lnab{Z}J_1(X),J_2X). 
\end{eqnarray*}
Assume that at a given point $\lnab{}X(p_0)=0$. Thus at $p_0$,
$g(\lnab{Z}J_1(X),J_2X)$  $=g(\lnab{Z}(J_1(X)),J_2X)$ $=
-g(J_1X,\lnab{Z}(J_2(X))$ $=-g(J_1X,\lnab{Z}J_2(X))$.
Hitchin (\cite{[Hi]}, lemma (6.8)) proved that
hyper-almost-K\"{a}hler structures are in fact hyper-K\"{a}hler.\qed
\section{The irreducible components of $R\in {\cal R }$}
\setcounter{Th}{0}
\setcounter{Pp}{0}
\setcounter{Cr} {0}
\setcounter{Lm} {0}
\setcounter{equation} {0}
\baselineskip .5cm
 ${\cal R}$ has an 
orthogonal decomposition
${\cal R}=\oplus_{i=1}^4{\cal R}_i$ into $O(V)$-invariants subspaces 
(\cite{[Si-T]}),
where 
\begin{eqnarray*}
{\cal R}_1 &=& (Ker~ b)^{\bot}=\{R:\sigma_R=0\}
=\{R: R=\lambda *\}\subset Ker~Ricci \subset Ker~tr\\[-1mm]
{\cal R}_2 &=& (Ker~ tr)^{\bot}=\{R=\lambda Id\}\\[-1mm]
{\cal R}_3 &=& (Ker~ b)\cap (Ker~ Ricci)=\{R: *R=R*, tr~R=tr* R=0\}\\[-1mm]
{\cal R}_4 &=& (Ker~ b)\cap (Ker~ Ricci)^{\bot}=\{R: *R=-R*\}
\end{eqnarray*}
and they can be characterized in terms of the sectional and Ricci curvature.
Set ${\cal B}:= Ker~b={\cal R}_2\oplus {\cal R}_3\oplus {\cal R}_4$.
We also have the following characterization using $K_{isot}$:
\begin{Pp} If $R\in {\cal R}$ and $b(R)=0$ then
$K_{isot}=0$ iff $\sigma_R(P)=-\sigma_R(P^{\bot})~\forall
P$. Hence, \\[-3mm]
\[{\cal R}_4=\{R: b(R)=0,~ \sigma_R(P)=-\sigma_R(P^{\bot})~\forall P\}
=\{R\in {\cal B}: K_{isot}=0\}=Sym_0(V)\bullet g\]
and $\phi\ra \phi\bullet g$ is an isomorphism from
 $Sym_0(V)$ onto  ${\cal R}_4$.
\end{Pp}
\em Proof. \em  In \cite{[Si-T]} it is proved
the first  equality of ${\cal R}_4$,  and that the elements
of this set satisfy  $R(P,P^{\bot})=0$.
Now it turns out that this condition and 
$\sigma_R(P)=-\sigma_R(P^{\bot})$ is equivalent to $K_{isot}=0$. 
We prove the less obvious implication.
$K_{isot}=0$ means that for any o.n.b.\ $X,Y,Z,W$,
$~2R(X,Y,Z,W)=R(X,Z,X,Z)+R(X,W,X,W)+R(Y,Z,Y,Z)+R(Y,W,Y,W)$.
Setting $P=span
\{X,Y\}$  
and replacing $X$ by $-X$ we conclude\\[-3mm]
\[\begin{array}{l}
\sm{0=R(X,Y,Z,W)= R(P,P^{\bot})}\\
\sm{0=R(X,Z,X,Z)+R(X,W,X,W)+R(Y,Z,Y,Z)+R(Y,W,Y,W)
=:\Psi(X,Y,Z,W).}~~~
\end{array}\]
From $\sm{0=\Psi(X,Y,Z,W)=\Psi(Z,Y,X,W)=\Psi(Z,X,W,Y)}$ we have
$\sm{0=\Psi(Z,Y,X,W)}$ $\sm{-\Psi(X,Y,Z,W)}$ $\sm{+\Psi(Z,X,W,Y)}$ $
\sm{=2R(Z,W,Z,W)}$ $\sm{+2R(X,Y,X,Y)}$ $
\sm{=2\sigma_R(P^{\bot})+2\sigma_R(P)}.$ For a proof that 
$\phi\ra \phi\bullet g$ defines an isomorphism see  e.g. \cite{[Si-T]}.
\qed\\

Let $\phi \in Sym(V)$, and $e_i$ a d.o.n.b. of eigenvectors of $\phi$, 
with corresponding  eigenvalues $\lambda_i$. Then
$e_i\wedge e_j$ with $i<j$ are the eigenvectors of
$\phi\bullet g$ ($\phi\bullet\phi$ resp.) corresponding to the respective 
eigenvalues $\lambda_i+\lambda_j$ ($2\lambda_i\lambda_j$ resp.).
Let $\Lambda^{\pm}_1=e_1\wedge e_2 \pm e_3\wedge e_4$,
$\Lambda^{\pm}_2=e_1\wedge e_3 \mp e_2\wedge e_4$,
$\Lambda^{\pm}_3=e_1\wedge e_4 \pm e_2\wedge e_3$.
We have:
\begin{eqnarray*}\begin{array}{cc}
\langle \phi\bullet g(\Lambda_{\al}^{\pm}), \Lambda_{\be}^{\pm}\rangle
= \delta_{\al\be}tr(\phi) &
\langle \phi\bullet g(\Lambda_{1}^+), \Lambda_{\be}^-\rangle
= \delta_{1\be}(\lambda_1+\lambda_2-\lambda_3-\lambda_4)\\
\langle \phi\bullet g(\Lambda_{2}^+), \Lambda_{\be}^-\rangle
= \delta_{2\be}(\lambda_1+\lambda_3-\lambda_2-\lambda_4)~~~~~~~~&
\langle \phi\bullet g(\Lambda_{3}^+), \Lambda_{\be}^-\rangle
= \delta_{3\be}(\lambda_1+\lambda_4-\lambda_2-\lambda_3)
\end{array}\\[-7mm]\end{eqnarray*}
\begin{eqnarray*}\begin{array}{cc}
\langle \phi\bullet\phi(\Lambda_{1}^{\pm}), \Lambda_{\be}^{\pm}\rangle
= 2\delta_{1\be}(\lambda_1\lambda_2+ \lambda_3\lambda_4)~~~~~~
&
\langle \phi\bullet\phi(\Lambda_{1}^{\pm}), \Lambda_{\be}^{\mp}\rangle
= 2\delta_{1\be}(\lambda_1\lambda_2- \lambda_3\lambda_4)\\
\langle \phi\bullet\phi(\Lambda_{2}^{\pm}), \Lambda_{\be}^{\pm}\rangle
= 2\delta_{2\be}(\lambda_1\lambda_3+\lambda_2\lambda_4)~~~~~~
&
\langle \phi\bullet\phi(\Lambda_{2}^{\pm}), \Lambda_{\be}^{\mp}\rangle
= 2\delta_{2\be}(\lambda_1\lambda_3-\lambda_2\lambda_4)\\
\langle \phi\bullet\phi(\Lambda_{3}^{\pm}), \Lambda_{\be}^{\pm}\rangle
= 2\delta_{3\be}(\lambda_1\lambda_4+\lambda_2\lambda_3)~~~~~~
&
\langle \phi\bullet\phi(\Lambda_{3}^{\pm}), \Lambda_{\be}^{\mp}\rangle
= 2\delta_{3\be}(\lambda_1\lambda_4-\lambda_2\lambda_3).\end{array}
\end{eqnarray*}
\begin{Lm}$
(\phi\bullet g)^{\pm}_{\pm}= \ha tr(\phi)Id_{\bigwedge^{\pm}_2},$ ~and~ 
 $tr(\phi\bullet \phi)=
\sum_{i<j}2\lambda_i\lambda_j=2\sigma_2(\phi).$
\end{Lm}
\begin {Lm} Let $\phi,\xi\in Sym(V)$ and $\langle \phi,\xi\rangle=
\sum_{ij}\phi(e_i,e_j)\xi(e_i,e_j)$ where $e_i$ is an o.n.b.\ of
$V$. Then\\[-8mm]
\begin{eqnarray*}\begin{array}{ccc}
\xi\bullet\phi \in {\cal B}&~~~~~~&
tr(\phi\bullet g)=\langle Id, \phi\bullet g\rangle = 3tr (\phi)\\
\langle\phi\bullet g,\xi\bullet g\rangle =2\langle \phi, \xi\rangle
+tr(\phi)tr(\xi)&~~~~~~&
*(\phi \bullet g) * =\ha tr(\phi) g\bullet g-\phi\bullet g\\
 Ricci_{\phi\bullet g}= Tr(\phi)g +2\phi &~~~~~~&
 s_{\phi\bullet g}=6tr \phi\\
\langle \phi \bullet g, *\xi \bullet g\rangle =0
&~~~~~~& \phi\bullet g=\xi\bullet g\mbox{~~iff~~} \phi=\xi.\end{array}
\end{eqnarray*}
\end{Lm}
\em Proof. \em  The equalities are proved using the eigenvalues
$\lambda_i$ and  eigenvectors $e_i$  of $\phi$. Since both
 $P=e_i\wedge e_j$ and $P^{\bot}$ are eigenvectors of $R=\phi\bullet g$,
then  $*R*$ has the same eigenvalues $\lambda'_i+\lambda'_j$
as $R$.  Namely, if $(ijkl)$ is a
permutation of $1234$,
$\lambda'_i+\lambda'_j= \sum_s\lambda_s -(\lambda_i+\lambda_j)
= tr(\phi)-(\lambda_i+\lambda_j)
=\lambda_k +\lambda_l$.\qed\\[4mm]
If $R\in {\cal R}$, then $R_1=\frac{1}{3}t(R)*=\frac{1}{6}tr(*R)*$.
For  $R\in{\cal B}$ we have the decomposition:
\begin{eqnarray}
R_2&=&\sm{\frac{1}{12}} tr(R) g\bullet g~=\sm{\frac{1}{24}} s_R\, 
 g\bullet g =\sm{\frac{1}{12}}s_R Id\\
R_3&=& \ha(R+*R*)-\sm{\frac{1}{12}} tr(R)g\bullet g~=~{\cal W}
={\cal W}^++{\cal W}^-\\
R_4&=& \ha(R-*R*)~=~
\ha (Ricci_R - \sm{\frac{1}{4}}s_R\, g)\bullet g.
\end{eqnarray}
 ${\cal W}$ is the Weyl tensor.  It  applies
$\bigwedge^2_{\pm}V$ into $\bigwedge^2_{\pm}V$, 
with ${\cal W}^{+}$ and ${\cal W}^{-}$ the self-dual  and
the anti-selfdual part respectively, i.e, they satisfy:
$~*{\cal W}^{\pm}={\cal W}^{\pm}*=\pm{\cal W}^{\pm}.$ 
Furthermore
$~R^{\pm}_{\pm}= {\cal W}^{\pm}+ \frac{s_R}{12}Id_{\bigwedge^{\pm}}$, ~
$R^+_-+ R^-_+ = \ha (Ricci - \frac{s_R}{4}g)\bullet g$. 
From previous lemmas we have:
\begin{Pp} ${\cal R}_4\oplus {\cal R}_2=Sym(V)\bullet g$, and $\forall \phi
\in Sym(V)$ 
$(\phi\bullet g)_2=\frac{tr(\phi)}{4}g\bullet g$, 
$(\phi\bullet g)_4=(\phi-\frac{tr(\phi)}{4}g)\bullet g$.
\end{Pp}
\begin{Lm} Let $R\in {\cal B}$ and  $\phi, \xi\in Sym(V)$. Then
~$\langle  R, \phi\bullet g\rangle = \langle Ricci_R,\phi\rangle$,~
$\langle  * R*, \phi\bullet g\rangle = \ha s_R\, 
tr(\phi)-\langle Ricci_R,\phi\rangle$, and 
~$\langle  * \phi\bullet g*, \phi\bullet g\rangle =  2(tr \phi)^2
-2\|\phi\|^2=2tr(\phi\bullet \phi).$
\end{Lm}
\em Proof. \em Using the decomposition (4.1)-(4.3), lemma 4.2
and proposition 4.2 we have
\begin{eqnarray*}
\langle R, \phi\bullet g\rangle &=& \sm{\frac{1}{24}}s_R\langle g\bullet g, 
\phi\bullet g\rangle + \ha \langle (Ricci_R-\sm{\frac{1}{4}}s_R\, g)\bullet g,
\phi\bullet g\rangle\\
&=&\sm{\frac{1}{24}}s_R(2\langle g,\phi\rangle + tr(g)tr(\phi))
+\langle (Ricci_R-\sm{\frac{1}{4}}s_R\, g), \phi\rangle
= \langle Ricci_R,\phi\rangle.\\[-5mm]
\end{eqnarray*}
\[
\langle *R*, \phi\bullet g\rangle = \langle Ricci_{*R*},\phi\rangle=
\langle \ha s_R\, g-Ricci_R, \phi \rangle
= \sm{\frac{1}{2}} s_R\,tr(\phi)-\langle Ricci_R, \phi\rangle.  \qed.\\[3mm]
\]
\begin{Pp} If $Q=R+ \phi\bullet g$, where $R \in {\cal B}$ and 
$\phi \in Sym(V)$, then
\begin{eqnarray}
4\pi^2{\cal X}(Q)&=& 4\pi^2{\cal X}(R)
+\LA{(}\frac{s_R}{2}tr(\phi)-\langle Ricci_R, \phi\rangle 
+(tr (\phi))^2-\|\phi\|^2 \LA{)}Vol_V \\
4\pi^2 p_1(Q)&=&4\pi^2 p_1(R).
\end{eqnarray}
\end{Pp}
\em Proof. \em From lemma 4.3\\[-6mm]
\begin{eqnarray*}
4\pi^2\langle{\cal X}(Q), Vol_V \rangle 
&=& \ha \langle *Q*, Q\rangle = \ha \langle *R*,R\rangle
+ \langle *R*, \phi\cdot g\rangle + \ha \langle *(\phi\bullet g)*,
\phi\bullet g\rangle \\
&=&4\pi^2\langle {\cal X}(R),Vol_V \rangle
+\sm{\frac{s_R}{2}}tr(\phi)-\langle Ricci_R, \phi\rangle
+(tr (\phi))^2 -\|\phi\|^2.
\end{eqnarray*}
Note that  $R=R^T$, $Q=Q^T$. 
Using decomposition (4.1)(4.2)(4.3) and lemmas 4.2, 4.1, we have
\begin{eqnarray*}
4\pi^2\langle p_1(Q),  Vol_V \rangle
&=&\langle Q,*Q\rangle =\langle R,*R\rangle
+2\langle R,*\phi\bullet g\rangle 
+\langle \phi\bullet g,*\phi\bullet g\rangle\\
&=&4 \pi^2\langle p_1(R), Vol_V \rangle
+2\langle {\cal W},*\phi\bullet g\rangle = 4\pi^2 p_1(R)
+2\langle *{\cal W},\phi\bullet g\rangle\\
&=&4 \pi^2\langle p_1(R), Vol_V \rangle 
+2\langle {\cal W}^+,\phi\bullet g\rangle
-2\langle {\cal W}^-,\phi\bullet g\rangle\\
&=&4 \pi^2\langle p_1(R), Vol_V \rangle
+2\langle {\cal W}^+,(\phi\bullet g)^+_+\rangle
-2\langle {\cal W}^-,(\phi\bullet g)^-_ -\rangle\\
&=&4 \pi^2\langle p_1(R), Vol_V \rangle 
+tr(\phi)\La{(}\langle {\cal W}^+,Id_{\Lambda^+_2}\rangle
-\langle {\cal W}^-,Id_{\Lambda^-_2}\rangle\La{)}\\
&=&4\pi^2 \langle p_1(R), Vol_V \rangle 
+tr(\phi)(tr({\cal W}^+)-tr({\cal W}^-)).
\end{eqnarray*}
Since ${\cal W}$ and $*{\cal W}$ have zero trace, the same holds for 
${\cal W}^+$ and ${\cal W}^-$ leading to (4.5). \qed\\[-2mm]
\begin{Cr} If $\lambda_i$ are the eigenvalues
of $\phi$ then ${\cal X}(\phi\bullet g)=0$ iff
$\sigma_2(\phi)\!:=\sum_{i<j}\lambda_i\lambda_j=0$.
\end{Cr}
\begin{Pp}
 If $R\in {\cal B}$ and $\lambda\in \R{}$, then\\[-2mm]
\[\begin{array}{l}
4\pi^2 {\cal X}(R+\lambda *) = 4\pi^2{\cal X}(R)+3\lambda^2Vol_V\\[1mm] 
4\pi^2 p_1(R+\lambda *)= 4\pi^2 p_1(R)$ $+  \lambda\, s_R Vol_V.
\end{array}\]
\end{Pp}
\em Remark $3$. \em Using (4.1)-(4.3) one sees that for $R\in {\cal B}$
(see e.g. \cite{[Be1]})\\[-3mm]
\[\begin{array}{l}
8\pi^2 {\cal X}(R)=(\|R_2\|^2+\|R_3\|^2-\|R_4\|^2)Vol_V\\
4\pi^2 p_1(R)=\langle R_3, R_3*\rangle Vol_V =(\|{\cal W}^+\|^2
-\|{\cal W}^-\|^2) Vol_V.
\end{array}\]
The invariance of the Pontrjagin class under a conformal change
of the metric now follows from the well known conformal invariance of
the Weyl tensor, but we give a simple proof in Prop.4.3.
\section{ Conformal change of the metric}
\setcounter{Th}{0}
\setcounter{Pp}{0}
\setcounter{Cr} {0}
\setcounter{Lm} {0}
\setcounter{equation} {0}
\baselineskip .5cm
 Let $g$ be a Riemannian metric on a 4-dimensional manifold $M$ and 
$f:M\ra \R{}$ a smooth  map. We denote by $\lnab{}$ and $\hlnab{}$ the
Levi Civita connections of $g$ and $\hat{g}=e^fg$ with
curvature tensors $R^M$, $\hat{R}^M$ respectively.
Then  $\hlnab{X}Y=\lnab{X}Y +\hat{S}(X,Y)$
where $\hat{S}(X,Y)=\ha\La{(} d f(X)Y +df(Y)X
-\langle X,Y\rangle \lnab{}f\La{)}$. 
The curvature tensors $\hat{R}^M(X,Y,Z,W)$ $=\hat{g}(\hat{R}^M(X,Y)Z,W)$,
and  ${R}^M(X,Y,Z,W)$ $={g}({R}^M(X,Y)Z,W)$ are related by
\begin{equation}
\hat{R}^M(X,Y,Z,W) =  e^f\, R^M(X,Y,Z,W)  +e^f\, \phi\bullet g(X,Y,Z,W)
\end{equation}
\begin{equation}
\mbox{where~~~~~~}\phi=\ha\La{(}-
\sm{\frac{1}{4}}\|\nabla f\|^2 g
+\ha df\otimes d f  -Hess (f)\La{)}.
\end{equation}
Let $Hess (f)(X,Y)= \lnab{X}df(Y)=\lnab{Y}df(X)$ be the Hessian of $f$,
 $\Delta df =(d\delta+\delta d)df= d\delta df$, 
and  $S(df)(Y)=df(Ricci^M(Y))=Ricci^M(Y,\nabla f)$, the
Weitzenb\"{o}ck operator $S$ on the 1-form $df$.
We are using the following sign for the Laplacian of maps
$h:M\ra \R{}$,~$\Delta h= \Delta^+ h=-\delta dh = div(\nabla h)
= tr(Hess(h))$. 
 Set
$\lnab{X,Y}^2df=\lnab{X}(\lnab{Y}df)-\lnab{\lnab{X}Y}df$.\\[-4mm]
\begin{Lm} For all $X,Y,Z\in T_pM$, $p\in M$ \\[3mm]
$(a)~~\lnab{X,Y}^2df(Z)=\lnab{X,Z}^2df(Y).$\\[1mm]
$(b)~~\lnab{X,Y}^2df(Z)-\lnab{Y,X}^2df(Z)=df(R^M(X,Y)Z)=
\lnab{X,Z}^2df(Y)-\lnab{Y,Z}^2df(X).$\\[1mm]
$(c)~~\Delta d f= -d(\Delta f).$\\[1mm]
$(d)~~-\ha\Delta\|\nabla f\|^2=\langle \Delta d f, d f\rangle
-\|Hess (f)\|^2 -Ricci^M(\nabla f, \nabla  f).$\\[1mm]
$(e)~~(\Delta f)^2=-\ha \Delta \|\nabla f\|^2+
div(\Delta f\, \nabla f)+\|Hess (f)\|^2+
Ricci^M(\nabla f,\nabla f).$\\[1mm]
$(f)~~\|\nabla f\|^2\, \Delta f=div(\|\nabla f\|^2\nabla f)
-2Hess (f) (\nabla f, \nabla f).$\\[-2mm]
\end{Lm}
\em Proof. \em  
We may assume without loss of generality
that $\lnab{}X=\lnab{}Y=\lnab{}Z=0$ at a given $p_0$.
From
$\lnab{X}\left(\lnab{Z}df(Y)\right)
=\lnab{X}\left(\lnab{Y}df(Z)\right)$, we get (a) at $p_0$.
At $p_0$\\[-5mm]
\begin{eqnarray*}
\lefteqn{\lnab{X}\left(\lnab{Y}df(Z)\right)-
\lnab{Y}\left(\lnab{X}df(Z)\right)=}\\
&=&\lnab{X}\left(\lnab{Y}(df(Z))\right)-
\lnab{Y}\left(\lnab{X}(df(Z))\right)-
\lnab{X}(df(\lnab{Y}Z))+\lnab{Y}(df(\lnab{X}Z))\\
&=&df(-\lnab{X}\lnab{Y}Z+\lnab{Y}\lnab{X}Z)=df(R^M(X,Y)Z).
\end{eqnarray*}
Thus (b) is proved.
By (b), if $e_1,\ldots, e_4$ is any o.n. frame of $M$
with $\lnab{}e_i(p_0)=0$, then at $p_0$
\[df(Ricci^M(Y))=\sum_i\lnab{e_i,e_i}^2df(Y)-\lnab{Y,e_i}^2df(e_i)=
\sum_i\lnab{e_i,e_i}^2df(Y)-\lnab{Y}(\Delta f).\\[-3mm]\]
Consequently, by the Weitzenb\"{o}ck formula\\[-3mm] 
\[
\Delta df(Y) =(d\delta+\delta d)df(Y)=
-(Trace\, \lnab{}^2df)(Y) +S(df)(Y)=-\lnab{Y}(\Delta f).
\]
So, we have proved (c). (d) is a well-known application of the
Weitzenb\"{o}ck formula for 1-forms.  (e) and (f) follows from
previous alignments and that $div(hX)=\langle X,\nabla h\rangle +h\, div(X)$.
\qed \\[1mm]
\begin{Lm}~\\[-9mm]
\begin{eqnarray}
tr(\phi)\! &=&\!\!-\sm{\frac{1}{4}}\|\nabla f\|^2
-\ha \Delta f\\[1mm]
tr(\phi\bullet \phi)=(tr(\phi))^2\!-\!\|\phi\|^2\!&=&
\!\sm{\frac{1}{8}}div\La{(}(2\Delta f
\!+\!\|\nabla f\|^2)\nabla f-\nabla(\|\nabla f\|^2)\La{)}
+\sm{\frac{1}{4}}Ricci^M(\nabla f, \nabla f)~~~~~~~~~\\[1mm]
\ha s^M tr(\phi)-\langle Ricci^M, \phi\rangle\! &=&
\!-\sm{\frac{1}{4}}s^M\Delta f-\sm{\frac{1}{4}}
 Ricci^M(\nabla f, \nabla f)
+\ha\langle Ricci^M, Hess (f)\rangle.\\[-2mm]\non
\end{eqnarray}
\end{Lm}
\em Proof. \em Let $e_i$ be an o.n.b. of $T_pM$. We have 
\\[-6mm]
\begin{eqnarray*}
tr(\phi) &=&\ha\La{(}-\|\nabla f\|^2+\sum_i\ha df(e_i)df(e_i)
-Hess(f)(e_i,e_i)\La{)}
=-\sm{\frac{1}{4}}\|\nabla f\|^2-\ha \Delta f\\[-6mm]
\end{eqnarray*}
\[
\|\phi\|^2=\sm{\frac{1}{4}}\|Hess f\|^2 +
\sm{\frac{1}{8}}\|\nabla f\|^2\Delta f
-\sm{\frac{1}{4}} Hess (f)(\nabla f, \nabla f)
+\sm{\frac{1}{16}\|\nabla f\|^4}.
\]
Thus\\[-5mm]
\[(tr(\phi))^2-\|\phi\|^2 =\sm{\frac{1}{4}}(\Delta f)^2
+\sm{\frac{1}{8}}\|\nabla f\|^2\Delta f-
\sm{\frac{1}{4}}\|Hess (f)\|^2
+\sm{\frac{1}{4}} Hess (f)(\nabla f, \nabla f).\]
Using Lemma 5.1 (e) and (f) we get (5.4).  Now
\begin{eqnarray*}
\lefteqn{\sm{\frac{1}{2}}s^M tr(\phi)-\langle Ricci^M, \phi\rangle =}\\
&=&\sm{\frac{1}{2}}s^M (-\sm{\frac{1}{4}}\|\nabla f\|^2
-\ha\Delta f)+ \sm{\frac{1}{8}}\|\nabla f\|^2s^M 
-\sm{\frac{1}{4}} Ricci^M(\nabla f, \nabla f)
+\sm{\frac{1}{2}}\langle Ricci^M, Hess (f)\rangle\\
&=&-\sm{\frac{1}{4}}s^M\Delta f
-\sm{\frac{1}{4}} Ricci^M(\nabla f, \nabla f)
+\ha\langle Ricci^M, Hess (f)\rangle.~~~~~~~~\qed\\[-1mm]
\end{eqnarray*}
We recall the  well-known
formula (see e.g.\ \cite{[G-H-L]}) $ \delta (Ricci^M)^{\sharp} 
+ \ha \nabla  s^M=0,$
derived from the second Bianchi
identity,
 where $\delta$ is the formal adjoint to $d$ for $TM$-valued forms on $M$. 
Now
\begin{eqnarray*}
div\La{(} (Ricci^M)^{\sharp}(\nabla f)\La{)} &=&
\sum_i \langle \lnab{e_i}((Ricci^M)^{\sharp}(\nabla f)), e_i\rangle\\[-1mm]
&=&\sum_i \langle \La{(}\lnab{e_i}(Ricci^M)^{\sharp}\La{)}(\nabla f), 
e_i\rangle+\langle  (Ricci^M)^{\sharp}(\lnab{e_i}\nabla f), e_i\rangle\\[-1mm]
&=&-\langle\delta((Ricci^M)^{\sharp}), \nabla f\rangle
+\langle Ricci^M, Hess(f)\rangle.\\[-4mm]
\end{eqnarray*}
Moreover, ~$g(\nabla  s^M, \nabla f) =\langle d s^M, df\rangle =d(s^M*df)
+s^M\delta df =div(s^M\nabla f)-s^M\Delta f$.
Thus
\begin{eqnarray*}
0 &=&g(\delta (Ricci^M)^{\sharp} + \ha \nabla  s^M, \nabla f)\\
&=& div \La{(} -(Ricci^M)^{\sharp}(\nabla f)+\ha s^M\nabla f\La{)}
+\langle Ricci^M, Hess(f)\rangle- \ha s^M\Delta f.
\end{eqnarray*}
So we have obtained:\\[-5mm]
\begin{Lm}
~~~$\langle Ricci^M, Hess(f)\rangle- \ha s^M\Delta f
=div(-\ha s^M\nabla f+ (Ricci^M)^{\sharp}(\nabla f))$.\\
\end{Lm}
\em Proof of Proposition 1.1.  \em 
For each vector field $Z$ on $M$ we denote
$\hat{Z}=e^{-\frac{f}{2}}Z$.  The
curvature tensor $\hat{R}^M$ of $TM$ w.r.t the Levi-Civita connection
of the Riemannian metric $\hat{g}$ can be seen as a
bundle map
$Q:\bigwedge^2TM\ra \bigwedge^2TM $  given by
$g(Q(X\wedge Y),(Z\wedge W))$ $=Q(X,Y,Z,W)$ $=
{\hat{R}^M}( X,Y, \hat{Z},\hat{W})$.
Thus $Q(X,Y,Z,W)= R^M(X,Y,Z,W)
+\phi\bullet g(X,Y,Z,W),$ 
where $\phi$ is given by (5.2). Therefore
${\cal X}(\hat{R}^M)={\cal X}(Q)={\cal X}(R^M+ \phi\bullet g).$ 
Now we apply Proposition 4.3  and Lemmas 5.2 and 5.3 
to obtain (1.2). Similar for $p_1(\hat{R}^M)$. \qed
\section{Connections with torsion}
\setcounter{Th}{0}
\setcounter{Pp}{0}
\setcounter{Cr} {0}
\setcounter{Lm} {0}
\setcounter{equation} {0}
\baselineskip .5cm
Let $(M,g)$ be a manifold of dimension 4 and $\lnab{}$ be its
Levi-Civita connection, and $\lnabp{}$ any  other connection.
Set $T(X,Y)=\lnabp{X}Y-\lnabp{Y}X-[X,Y]$, and $S(X,Y)=\lnabp{X}Y-\lnab{X}Y$.
Then $T(X,Y)=S(X,Y)-S(Y,X).$ From $\lnabp{X}g(Y,Z)=g(S(X,Y),Z)+g(S(X,Z),Y)$,
we see that  $\lnabp{}$ is a $g$-Riemannian connection iff
$g(S(X,Y),Z)=-g(S(X,Z),Y)$.  Moreover
 $\lnabp{}= \lnab{}$ iff $T=0$, that is $S(X,Y)=S(Y,X)$.
If $\zeta:TM\ra E$ is a vector bundle map and $E\ra M$ has a connection 
$\lnabe{}$ then $d\zeta= d^{\lnabe{}}\!\!\!\zeta$ is a 2-form on 
$M$ with values on $E$,
$d\zeta(X,Y)=\lnabe{X}(\zeta(Y))-\lnabe{Y}(\zeta(X))
-\zeta([X,Y]). $
Since $\lnab{}$ is a torsion free connection on $M$, 
$d\zeta(X,Y)=\lnab{X}\zeta(Y)-\lnab{Y}\zeta(X).$
We define
${\cal S}:TM\ra L(TM;TM)$ and
${\cal S}^2:\bigwedge^2 TM\ra L(TM;TM)$  by
\[{\cal S}(X)(Y)=S(X,Y)~~~~~~~
{\cal S}^2(X\wedge Y)(Z)=S(X,S(Y,Z))-S(Y,S(X,Z)) \]
and take the exterior derivative
 $d{\cal S}:\bigwedge^2TM \ra L(TM;TM)$, by taking the usual connection
on $L(TM;TM)$ w.r.t the Levi-Civita connection of $M$. We have
\[
R'(X\wedge Y)(Z) = R^M(X\wedge Y)(Z)
-d{\cal S}(X\wedge Y)(Z) - {\cal S}^2(X\wedge Y)(Z).
\]
If $\lnabp{}\,$ is $g$-Riemannian then ${\cal S}$, $d{\cal S}$ and ${\cal S}^2$
take values on $\bigwedge^2TM$ and 
$\langle {\cal S}(X), Y\wedge Z\rangle = g(S(X,Y),Z)$,
\begin{eqnarray*}
\sm{\langle{\cal S}^2(X\wedge Y), Z\wedge W\rangle} &=&
\sm{g( S(X,Z),S(Y,W))-g(S(X,W),S(Y,Z))}\\
\sm{\langle R'(X\wedge Y),Z\wedge W\rangle }&=& 
\sm{\langle R^M(X\wedge Y), Z\wedge W\rangle
-\langle d{\cal S}(X\wedge Y),Z\wedge W\rangle 
-\langle {\cal S}^2(X\wedge Y), Z\wedge W\rangle.}
\end{eqnarray*}
\section{Transgression forms}
\setcounter{Th}{0}
\setcounter{Pp}{0}
\setcounter{Cr} {0}
\setcounter{Lm} {0}
\setcounter{equation} {0}
\baselineskip .5cm
Let $(M,g)$ be a Riemannian manifold of dimension 4 with
Levi-Civita connection $\lnab{}$ and let $\lnabp{}$ be any other
 $g$-Riemannian connection. 
We take $\pi: \tilde{M}=M\times [0,1]\ra M$ the projection map 
$\pi (p,t)=p$ and the connection  $\tlnab{}$ on $\pi^{-1}TM$ given by
$\tlnab{}= t\lnab{}+(1-t)\lnabp{}$, 
where $\lnab{}$ and $\lnabp{}$ are the pullback connections w.r.t $\pi$.
For any vector fields $X,Y,Z,W$ of $M$, map  $f:M\times [0,1]\ra \R{} $ 
and $u\in T_pM$, $h\in \R{}$,
$\tlnab{(u,h)}(f\pi^{-1}Y)_{(p,t)}$ $ 
=df(p,t)(u,h)Y_p$$ +f(p,t)\lnab{u}Y_{(p)}$$+tf(p,t)S_p(u,Y_p)$.
Thus
\begin{eqnarray*}
\tlnab{\ddt}\tlnab{(X,0)} (\pi^{-1}Y) &=& \tlnab{\ddt}\La{(}
\pi^{-1}(\lnab{X}Y) + t \pi^{-1}(S(X,Y))\La{)}=S(X,Y)
\end{eqnarray*}
and assuming at a given point $p$, $\lnab{}Y_p=\lnab{}Z_p=0$, then at $(p,t)$
\begin{eqnarray*}
\tlnab{(X,0)}\tlnab{(Y,0)} (\pi^{-1}Z)
&=& \lnab{X}\lnab{Y}Z + t\lnab{X}S(Y,Z) + t^2S(X,S(Y,Z)).
\end{eqnarray*}
Therefore, the curvature tensor of $\pi^{-1}TM$
with respect to this connection (Riemannian w.r.t $\pi^{-1}g$),  
$\tilde{\underline{R}}:\bigwedge^2\pi^{-1}TM \ra \bigwedge^2T\tilde{M}$~~
$\tilde{\underline{R}}(\pi^{-1}Z,\pi^{-1}W)(\tilde{X},\tilde{Y})
=:\tilde{R}(\tilde{X},\tilde{Y},\pi^{-1}Z,\pi^{-1}W)$
satisfies
\begin{eqnarray*}
\tilde{R}(X,Y,Z,W)&=&
R^M(X,Y,Z,W)-t d{\cal S}(X\wedge Y)(Z\wedge W)
-t^2{\cal S}^2(X\wedge Y)(Z\wedge W)\\
\tilde{R}(\ddt,Y,Z,W)&=&-{\cal S}(Y)(Z\wedge W)
\end{eqnarray*}
with $d{\cal S}=d^{\lnab{}}{\cal S}$ where $\lnab{}$ denotes
de induced connection on $\bigwedge^2TM$ by $\lnab{}$.
Let $e_1,e_2, e_3, e_4$ be an d.o.n.b.\ of $T_pM=(\pi^{-1}TM)_{(p,t)}$.
The definitions of ${\cal X}(\tilde{R})$ and $\pi_1(\tilde{R})$ are 
the 4-forms on $\tilde{M}$
\[
{\cal X}(\tilde{R}) 
=\sum_{i<j}\frac{1}{8\pi^2}\underline{\tilde{R}}(e_i\wedge e_j)\wedge 
\underline{\tilde{R}}(*e_i\wedge e_j)~~~~~
p_1(\tilde{R})
=\sum_{i<j}\frac{1}{4\pi^2}\underline{\tilde{R}}(e_i\wedge e_j)\wedge
\underline{\tilde{R}}(e_i\wedge e_j).
\]
For $\zeta\!:\!V\!\ra\! \bigwedge^2V$ and $R\!\in 
\!L(\bigwedge^2V;\bigwedge^2V)$,we define a 3-form
\begin{equation}
\langle \zeta\wedge R\rangle (X,Y,Z):=\langle \zeta(X),R(Y,Z)\rangle
+\langle \zeta(Z),R(X,Y)\rangle + \langle \zeta(Y),R(Z,X)\rangle.
\end{equation}
Now we compute ${\cal X}(\tilde{R})(\ddt,X,Y,Z)$.
\begin{eqnarray*}
\lefteqn{8\pi^2{\cal X}(\tilde{R})(\ddt,X,Y,Z)=}\\
=\sum_{i<j}
&\LA{(}&\langle\tilde{R}(\ddt,X),e_i\wedge e_j\rangle
\langle\tilde{R}(Y,Z),* e_i\wedge e_j\rangle
-\langle\tilde{R}(\ddt,Y),e_i\wedge e_j\rangle
\langle\tilde{R}(X,Z),* e_i\wedge e_j\rangle\\[-3mm]
&&+\langle\tilde{R}(\ddt,Z),e_i\wedge e_j\rangle
\langle\tilde{R}(X,Y),* e_i\wedge e_j\rangle
+\langle\tilde{R}(Y,Z),e_i\wedge e_j\rangle
\langle\tilde{R}(\ddt,X),* e_i\wedge e_j\rangle\\[-1mm]
&&-\langle\tilde{R}(X,Z),e_i\wedge e_j\rangle
\langle\tilde{R}(\ddt,Y),* e_i\wedge e_j\rangle
+\langle\tilde{R}(X,Y),e_i\wedge e_j\rangle
\langle\tilde{R}(\ddt,Z,* e_i\wedge e_j\rangle
~~~\LA{)}\\[2mm]
=\sum_{i<j}&2\LA{(}&\sm{-{\cal S}(X)(e_i\wedge e_j)
\La{(}\langle R^M(Y,Z),*e_i\wedge e_j\rangle
-t\langle d{\cal S}(Y\wedge Z),*e_i\wedge e_j\rangle
- t^2\langle{\cal S}^2(Y\wedge Z),*e_i\wedge e_j\rangle}\La{)}\\[-4mm]
&&\sm{+{\cal S}(Y)(e_i\wedge e_j)
\La{(}\langle R^M(X,Z),*e_i\wedge e_j\rangle
-t\langle d{\cal S}(X\wedge Z),*e_i\wedge e_j\rangle
- t^2\langle {\cal S}^2(X\wedge Z, *e_i\wedge e_j\rangle}\La{)}\\[-1mm]
&&\sm{-{\cal S}(Z)(e_i\wedge e_j)\La{(}\langle R^M(X,Y),*e_i\wedge e_j\rangle
-t\langle d{\cal S}(X\wedge Y),*e_i\wedge e_j\rangle
- t^2\langle {\cal S}^2(X\wedge Y),*e_i\wedge e_j\rangle\La{)}}~~\LA{)}~~~\\
&=&\sm{-2\langle {\cal S}\wedge *\La{(}R^M-td{\cal S}
- t^2{\cal S}^2\La{)}\rangle (X,Y,Z).}
\end{eqnarray*}
Similarly, $4\pi^2 p_1(\tilde{R})(\ddt,X,Y,Z)=\sm{-2\langle {\cal S}
\wedge \La{(}R^M-td{\cal S}- t^2{\cal S}^2\La{)}\rangle (X,Y,Z)}$.
Integration  over $[0,1]$ gives:
\begin{Pp} If $(M,g)$ is an oriented Riemannian manifold of dimension 4 with
Levi-Civita connection $\lnab{}$, and
$\lnabp{X}Y=\lnab{X}Y +S(X,Y)$ is a $g$-Riemannian connection on $TM$
with torsion $T(X,Y)=S(X,Y)-S(Y,X)$ and curvature tensor $R'$,
then
\begin{eqnarray}
{\cal X}(R') &=& {\cal X}(R^M)- \frac{1}{4\pi^2}d\La{(}\langle {\cal S}\wedge
*\la{(}R^M -\frac{1}{2}d{\cal S}- \frac{1}{3}{\cal S}^2\la{)}\rangle\La{)}\\
p_1(R')&=&  p_1(R^M) -\frac{1}{2\pi^2}d\La{(}\langle {\cal S}\wedge
\la{(}R^M -\frac{1}{2}d{\cal S}
- \frac{1}{3}{\cal S}^2\la{)}\rangle\La{)}.
\end{eqnarray}
\end{Pp}
\em Remark $4$. \em The Chern-Simons transgression formula is a particular 
case of the following fact: If $\varphi$ is a closed $k$-form on 
$M\times [a,b]$ and $\varphi_t= i_t^*\varphi$, where 
$i_t:M\ra M\times [a,b]$, $i_t(p)=(p,t)$, then 
$\varphi_b-\varphi_a= d T$
where $T$ is the $(k-1)$-form on $M$ given by:
$ T(X_1,\ldots,X_{k-1})=
\int^b_a \varphi(\ddt,(X_1,0), \ldots,
(X_{k-1},0))dt.$
\section{Singular connections}
\setcounter{Th}{0}
\setcounter{Pp}{0}
\setcounter{Cr} {0}
\setcounter{Lm} {0}
\setcounter{equation} {0}
\baselineskip .5cm
 Let $(M,g)$ be a Riemannian 4-manifold with Levi-Civita connection
$\lnab{}$ and $(E,g_E,\lnabe{})$ a Riemannian vector bundle
of rank 4 and Riemannian connection $\lnabe{}$ and assume there
exist a vector bundle map $\Phi:TM\ra E$ that is conformal
with $g_E(\Phi(X),\Phi(Y))=h(p)g(X,Y)$ $\forall X,Y\in T_pM$, 
where $h:M\ra \R{}$ is a smooth non-negative function with
zero set $\Sigma$. On $TM$ we have a degenerated metric $\hat{g}=
h g$, with Levi-Civita connection $\hlnab{}$, defined
on $M\sim \Sigma$ and let $\lnabp{}$ be the $\hat{g}$-Riemannian
connection on $TM$, defined along $M\sim \Sigma$, and that makes
$\Phi:TM \ra E$ a parallel isometry. We define
$\hat{S}(X,Y)= \hlnab{X}Y-\lnab{X}Y$,
$S'(X,Y)=\lnabp{X}Y-\hlnab{X}Y$, and
$S(X,Y)=$ $ \lnabp{X}Y-\lnab{X}Y$ $=S'(X,Y)+\hat{S}(X,Y).$
The connection $\lnabp{}$ satisfies $\Phi(\lnabp{X}Y)=\lnabe{X}(\Phi(Y))$,
and may not be smoothly extended to $\Sigma$, defining
a singular connection. But
the curvature tensor $R':\bigwedge^2TM\ra \bigwedge^2TM$ satisfies
\begin{eqnarray*}
R'(X,Y,Z,W)&=&\hat{g}(R'(X,Y)Z,W)= hg(R'(X,Y)Z,W)=
g_E(\Phi(R'(X,Y)Z),\Phi(W))\\
&=&g_E(R^E(X,Y)(\Phi(Z)), \Phi(W))= R^E(X,Y,\Phi(Z),\Phi(W))
\end{eqnarray*}
where $R^E:\bigwedge^2TM\ra \bigwedge^2E$ is the curvature tensor
of $(E,g_E, \lnabe{})$.
The above equality means that we can smoothly extend  $R'$ as a 4-tensor
to all $M$, defining it to be zero along $\Sigma$. Moreover, we
can define ${\cal X}(R')$ as $ {\cal X}(R^E)$ on
all $M$, and similar for $p_1(R')$.
Note that the star operator $*:\bigwedge^2TM\ra \bigwedge^2TM$
is the same for any conformal change of the metric.
Now we have
$\Phi(S(X,Y))=\Phi(\lnabp{X}Y)-\Phi(\lnab{X}Y)=
\lnabe{X}(\Phi(Y))-\Phi(\lnab{X}Y)=\lnab{X}\Phi(Y).$
Therefore
\begin{eqnarray}
S=\Phi^{-1}\circ \lnab{}\Phi~~~~~~~~~~~~T=\Phi^{-1}\circ d\Phi
\\
S'=S-\hat{S}= \Phi^{-1}\circ \lnab{}\Phi-\ha d\log h\odot Id
+\ha  g\otimes \nabla \log h
\end{eqnarray}
and the torsion $T'=T$.
From the two previous sections we obtain  the formulae in Theorem 1.1.
\begin{Pp} In the condition of Theorem 1.1, if $d\Phi=0$, then  
${\cal X}(R^E)={\cal X}(R^M)+\frac{1}{32\pi^2}div_g(P(\nabla \log h))Vol_M$ 
and $p_1(R^E)=p_1(R^M)$.
Moreover, in this case,  $\delta \Phi=0$ iff $h$ is constant, and so
${\cal X}(R^E)={\cal X}(R^M)$. This is the case when $\Phi$ is harmonic
and $M$ is compact.
\end{Pp}
\em Proof. \em 
If $d\Phi=0$, from (8.1), $T'=T=0$, and so
$\lnabp{}=\hlnab{}.$ Hence, $S'=0$, that is
$\Phi^{-1}\lnab{X}\Phi (Y)$ $=\ha d\log h(X)Y$ $+\ha d\log h (Y)X$ $-\ha
g(X,Y)\nabla \log h$. Thus,
$\Phi^{-1}(\delta \Phi)=\nabla \log h$, and so $\delta \Phi=0$
iff $h$ is constant. In this case ${\cal X}(R^E)={\cal X}(R^M)$,
$p_1(R^E)=p_1(R^M)$.
If $M$ is compact and $\Phi$ is harmonic,
 $\Delta \Phi=0$, then $d\Phi=0$ and $\delta \Phi=0$. \qed
\section{Homotopic almost complex structures}
\setcounter{Th}{0}
\setcounter{Pp}{0}
\setcounter{Cr} {0}
\setcounter{Lm} {0}
\setcounter{equation} {0}
\baselineskip .5cm
In this section we consider a connected compact oriented 
Riemannian 4-manifold $(M,g)$ with its Levi-Civita connection 
and we investigate  obstructions for two $g$-orthogonal 
almost complex structures
to be homotopic. This has been studied in \cite{[C]}, in a different context,
 using  shadows.
 We will use our framework. 
If $J_0$ and $J_1$ are
two positive $g$-orthogonal almost complex structures
\begin{equation}
J_1= \cos\theta(p)J_0 +\tilde{H}(p)=\cos\theta(p)J_0+ \sin\theta(p)H(p)
\end{equation}
where $\theta(p)\in[0,\pi]$ is  the \em angle \em
 function between $J_0$ and $J_1$, uniquely determined by  the
smooth function on $M$~  $\cos\theta(p)=
\ha\langle \omega_{J_0},\omega_{J_1}\rangle$, and
$\tilde{H}(p)$ is a smooth section on $E_{J_0}$ defined on all $M$.
Its normalized section $H$ may not be defined at the zero set 
${\cal Z}=\Sigma\cup
\Gamma$ of
$\tilde{H}$, where $\Sigma= J_1\cap -J_0=\{p\in M: J_1(p)=-J_0(p)\} $
is the set of \em anti-complex \em points
and $\Gamma=J_1\cap J_0=\{p\in M: J_1(p)=J_0(p)\}$ is the set of
\em common complex \em points.
 Any existent homotopy from $J_0$ to $J_1$ along 
$g$-orthogonal almost complex structures must be of the form
$J_t(p)=\cos\theta(p,t)J_0(p) +\tilde{H}(p,t)$.
If $\Sigma=\emptyset$ then $\tilde{H}(p)$ only vanish at points $p$ with
$\theta(p)=0$. In this case $J_1$ is homotopic to $J_0$.
For example,  we may take 
\begin{eqnarray}
\tilde{H}(p,t)&=& \sm{\sqrt{\frac{
t\la{(}2-t(1-\cos\theta(p))\la{)}}{(1+\cos\theta(p))}}}\tilde{H}(p)
=:\sin\theta(p,t) H(p) 
=\frac{\sin\theta(p,t)}{\sin\theta(p)}\tilde{H}(p)\\
\cos\theta(p,t)&=& 1-t(1-\cos\theta(p))\\
J_t &=&\cos\theta(p,t)J_0 +\tilde{H}(p,t)=
\cos\theta(p,t)J_0 +\sin\theta(p,t){H}(p).
\end{eqnarray}
Note that in this case 
$\cos\theta(p,t)$ is smooth on $M\times [0,1]$, $\tilde{H}$  smooth 
on $M\times ]0,1]$ and continuous on $M\times [0,1]$, but
not necessarily smooth, since $\frac{\sin\theta(p,t)}{\sin\theta(p)}$ 
may only  be continuous  at $t=0$. 
If we reparameterize $t$ on  (9.2) and (9.3)
  by a smooth function $\tau(t)$ 
with  $\tau(t)=e^{-\frac{2}{t}}$, $\tau(0)=0$
on a neighbourhood of $t=0$, then $J_t$ becomes smooth at $t=0$.
Note  also that $\sin\theta(p,t)
=\frac{\|\tilde{H}\|}{\sqrt{2}}$ (norm in $\bigwedge^2_+$) may 
only be continuous for $p\in\Gamma\cup \Sigma$. 
Hence, we have:
\begin{Pp} Any two positive  $\, g$-orthogonal almost complex structures, 
everywhere linearly independent, except  at common complex 
points, are homotopic.
\end{Pp}
Nevertheless, for generic $J_1$ such set $\Sigma$ is a
smooth surface. In fact (see \cite{[C]}) any almost complex structure $J$
 can be seen as an embedded  4-submanifold $J(M)$
 into the 6-dimensional  manifold $U\!\bigwedge^2_+TM$, the
total space of the sphere bundle (of radius $\sqrt{2}$) of $\bigwedge^2_+TM$,
 by  $p\ra (p,\omega_J(p))$. Since $\Sigma$ is the intersection
set of $J_1(M)$ with $-J_0(M)$, this is a surface for generic
$J_1$. Moreover, an orientation, depending on $J_1$ and $J_0$,
 can be given on each connected component of
$\Sigma$ (\cite{[C]}). \\[2mm]
$\Sigma$ can be seen as an obstruction to the existence of
a homotopy from $J_0$ to $J_1$, but we can always consider 
a smooth family of $g$-orthogonal almost complex
structures $J_t$, defined away from $\Sigma$ as in (9.4).
Then we consider an almost complex structure $\tilde{J}$
on ${\pi}^{-1}T(M\sim \Sigma)$, where $\pi:M\times [0,1]\ra M$,
$\pi(p,t)=p$, 
defined by $\tilde{J}_{(p,t)}={(J_t)}_p$.
Let $\tilde{\omega}_{(p,t)}=(\omega_{J_t})_p=\omega_{\tilde{J}}$,
 and define 
the bundle $\tilde{E}=E_{\tilde{J}}$ by 
$\tilde{E}_{(p,t)}=(E_{J_t})_p=:(E_t)_p$. Then
 $\pi^{-1}\bigwedge^2_+TM= \R{}\tilde{\omega}\oplus \tilde{E}$.
For $X\!\in\! T_pM$
we have $\lnab{(X,0)}^{\pi^{-1}}\omega_{\tilde{J}}(p,t)=
\lnab{X}^{}\omega_{{J}_t}$, where $\lnab{}^{\pi^{-1}}$ is the pullback
connection  on ${\pi}^{-1}\bigwedge^2_+TM$, and for $f$ a
function on $M\times [0,1]$ and $s$ a section of $\bigwedge^2_+TM$,~
$\lnab{(X,0)}^{{\tilde{E}}}(f \pi^{-1}s)(p,t)
=\lnab{X}^{E_t}(f(\cdot, t)s)$. 
Moreover $R^+(\ddt,(X,0))
=R^{M\times [0,1]}(\ddt,(X,0))=0$,
where now $\lnab{}^+$ denotes the induced
connection on $\pi^{-1}\bigwedge^+_2TM$ and  $R^+$ is its curvature
tensor.
Similar equations to (3.8) and (3.9) hold replacing $J$ by
$\tilde{J}$, and one has for $X,Y\in T_pM$,
 $ \eta_{\tilde{J}}(p,t)(X,Y)=
\eta_{J_t}(p)(X,Y)$.
Therefore
\begin{eqnarray}
c_1(\tilde{E})_{(p,t)}((X,0),(Y,0)) &=& c_1(E_t)_p(X,Y),\\[-1mm]
c_1(\tilde{E})_{(p,t)}(\ddt,(X,0)) &=&
\sm{\frac{1}{2\pi}}\eta_{\tilde{J}}(p,t)(\ddt, (X,0)).
\end{eqnarray}
So we have a closed 2-form $c_1(\tilde{E})$ on $(M\sim \Sigma)\times [0,1]$ 
with $i_t^*c_1(\tilde{E})=c_1(E_t)$ where $i_t(p)=(p,t)$. From Prop.\ 
3.1 and remark 4 of  section 7 we have proved  a first conclusion:
\begin{Pp}
If $J_0$ and $J_1$ are two positive homotopic $g$-orthogonal
almost complex structures on $M$, then as real cohomology classes
 $c_1(M,J_0)=c_1(M,J_1)$.  In particular, if $H^2(M;Z)$ has no torsion (
for example $M$ is simply connected) then
 $E_{J_0}$ and $E_{J_1}$ are isomorphic complex line bundles over $M$.
\end{Pp}
At $p\notin \Sigma\cup \Gamma$ we define, $\omega_2(0)_p=H(p)$, 
$\omega_3(0)_p=J_0(p)H(p)$, $\omega_1(0)_p=\omega_0$. Then $\omega_1(t)_p$
$:=\tilde{\omega}(p,t)=\cos\theta(p,t)\omega_1(0)_p+
\sin\theta(p,t)\omega_2(0)_p$, and
 a d.o.n.b.\ (of norm $\sqrt{2}$)
of $\tilde{E}_{(p,t)}$ is given by
\[
\omega_2(t)_p=-\sin\theta(p,t)\omega_1(0)_p
+\cos\theta(p,t)\omega_2(0)_p,~~~~~~~~~
\omega_3(t)_p=\omega_3(0)_p.
\]
Thus,\\[-8mm]
\begin{eqnarray*}
\lnab{(X,0)}^+\tilde{\omega}(p,t)&=& 
\cos\theta(p,t)\lnab{X}^+\omega_1(0)_p
+\sin\theta(p,t)\lnab{X}^+\omega_2(0)_p\\
&&+d \cos\theta(p,t)(X,0)\omega_1(0)_p+
d\sin\theta(p,t)(X,0)\omega_2(0)_p\\
\lnab{\ddt}^+\tilde{\omega}(p,t)&=&
{\ddt} \cos\theta(p,t)\omega_1(0)_p
+{\ddt}\sin\theta(p,t)\omega_2(0)_p\\[-5mm]
\end{eqnarray*}
\begin{eqnarray*}
\lefteqn{\eta_{\tilde{J}}({\ddt}, (X,0))
=\ha Vol_{\tilde{E}}\LA{(}\lnab{\ddt}^+
\tilde{\omega},\lnab{(X,0)}^+\tilde{\omega}(p,t)\LA{)}=}\\
&&~~~~~=~\LA{(}-\sin\theta(p,t)\ddt \cos\theta(p,t)
+\cos\theta(p,t)\ddt\sin\theta(p,t)\LA{)}\\ 
&&~~~~~~~~~~~~\cdot\LA{(}\cos\theta(p,t)
\langle \lnab{X}J_0, J_0H(p)\rangle
+\sin\theta(p,t)\langle \lnab{X}H, J_0H(p)\rangle\LA{)}
~~~~~~~~~~~~~~~~~~~
\end{eqnarray*}
where $\lnab{X}J$ is the covariant derivative in $Skew(TM)$ and
$\langle,\rangle$ is half the usual Riemannian metric on that vector bundle
(that corresponds to the Riemannian metric of $\bigwedge^2TM$).
For simplicity from now on we will denote by $\lnab{}$ the 
connection $\lnab{}^+$ of $\bigwedge^2_+TM$.
Note that
\begin{eqnarray*}
\lefteqn{
-\cos\theta\sin\theta \ddt\cos\theta
+\cos^2\theta\ddt\sin\theta=
-\ha\sin\theta\ddt\cos^2\theta
+\cos^2\theta\ddt\sin\theta}\\
&=& \ha \sin\theta\ddt \sin^2\theta +
\cos^2\theta \ddt\sin\theta
=\sin^2\theta\ddt\sin\theta+
\cos^2\theta \ddt \sin\theta
= \ddt\sin\theta.
\end{eqnarray*}
Similarly
$\cos\theta\sin\theta \ddt\sin\theta
-\sin^2\theta \ddt\cos\theta 
=-\ddt \cos\theta.$~
Consequently
\begin{eqnarray}
c_1(\tilde{E})_{(p,t)}(\ddt,(X,0)) &=&
\frac{1}{2\pi}\LA{(}-\ddt \cos\theta\, \langle
\lnab{X}H,J_0H\rangle +\ddt \sin\theta\, \langle 
\lnab{X}J_0,J_0H\rangle\LA{)}
\end{eqnarray}
and we have got the following formula:\\[-3mm]
\begin{Pp} If $J_0$  and $J_1$ are two  $g$-orthogonal almost complex
structures on $M$, with $J_1=\cos\theta J_0+\tilde{H}$
where $\tilde{H}\in C^{\infty}(E_{J_0})$, then
\begin{equation}
4\pi \La{(}c_1(M,J_1)-c_1(M,J_0)\La{)}= d\tilde{T} + dG
\end{equation}
where $\tilde{T}$ and $G$ are 1-forms on $M$, $G$ globally defined on $M$ and
$\tilde{T}$ defined away from $\Sigma$ by
\begin{equation}
\forall X\in T_pM~~~~~
~\tilde{T}(X)=\frac{1}{(1+\cos\theta(p))}\langle \lnab{X}\tilde{H},
J_0\tilde{H}\rangle, ~~~~~~
G(X)= \langle \lnab{X}J_0,J_0\tilde{H}\rangle.
\end{equation}
where now $\langle, \rangle$ is the usual inner product of
$Skew(TM)\subset TM^*\otimes TM$ (twice the one of $\bigwedge^2_+TM$).
\\[-4mm]
\end{Pp}
\em Proof. \em  We apply  remark 4
of section 7 to $c_1(\tilde{E})$ given by (9.5), (9.6), and satisfying
(9.7),  smooth on $(M\sim \Sigma\cup \Gamma)\times [0, 1]$,
obtaining $T(X)$ defined for $p\notin \Sigma\cup \Gamma$,
\begin{eqnarray*}
T(X) &=& 4\pi \int_{0}^1\!\!c_1(\tilde{E})
(\ddt,(X,0))dt = \int_{0}^1\!\! -\ddt 
\cos\theta(p,t)\langle \lnab{X}H,J_0H\rangle + 
\ddt\sin\theta(p,t)
\langle \lnab{X}J_0,J_0H\rangle~~~~~~~~~~\\[-1mm]
&=& (1-\cos\theta(p))\langle \lnab{X}H,J_0H\rangle
+\sin\theta(p)\langle \lnab{X}J_0,J_0H\rangle
=\tilde{T}(X)+ G(X).
\end{eqnarray*}
$T$ can be smoothly extended to $ \Gamma$.
 $\tilde{T}$ may not be defined
at $\Sigma$, but $G$ is smooth on all $M$. \qed
\begin{Cr} 
\begin{eqnarray*}
\\[-9mm]
\int_{M}c_1(M,J_1)\wedge c_1(M,J_0) &=& \int_M p_1(\sm{\bigwedge^2_+}TM)
-\frac{1}{32 \pi^{2}}\int_M d\La{(}\tilde{T}\wedge d (\tilde{T}+2G)\La{)}.
\end{eqnarray*}
\end{Cr}
\em Proof. \em
From Prop.3.1. $16 \pi^2 \la{(}c_1(M,J_1)-c_1(M,J_0)\la{)}^2 $ $=
16 \pi^2 \La{(} 2 p_1(\bigwedge^2_+ TM)-2 c_1(M,J_1)\wedge c_1(M,J_0)
\La{)}$.
Now formula of Prop. 9.3  gives \\[-4mm]
\[c_1(M,J_1)\wedge c_1(M,J_0) =p_1(\sm{\bigwedge^2_+}TM)
-\frac{1}{32 \pi^{2}}d\La{(}\tilde{T}\wedge d \tilde{T}
+ {G}\wedge d {G}+ 2\tilde{T}\wedge d G\La{)}. \qed \\[4mm]\]
\em Proof of Theorem 1.2. \em  Recall that
 $\omega_0$ is a closed 2-form. From  Prop.9.3
\begin{eqnarray*}
4\pi(c_1(M,J_1)-c_1(M,J_0))\wedge\omega_0 &=& d\tilde{T}\wedge \omega_0
+dG\wedge \omega_0 
= d(\tilde{T}\wedge \omega_0) + d(G \wedge \omega_0).
\end{eqnarray*}
Since $\tilde{T}\wedge \omega_0$ is a 3-form, 
$d(\tilde{T}\wedge \omega_0)
=-div\LA{(}(\star (\tilde{T}\wedge \omega_0))^{\sharp}\LA{)}Vol_M,$
where $\star:\bigwedge^3TM^*\ra TM^*$ is the star operator.
Now $\omega_0=e_*^{12}+e_*^{34}$  where $e_2=J_0e_1$, $e_4=J_0e_3$, and so
\begin{eqnarray*}
\star (\tilde{T}\wedge \omega_0) 
&=&\star \LA{(}\tilde{T}(e_3)e_*^{312}+\tilde{T}(e_4)e_*^{412}
+\tilde{T}(e_1)e_*^{134}+\tilde{T}(e_2)e^{234}\LA{)}\\
&=& -\tilde{T}(J_0e_4)e_*^4-\tilde{T}(J_0e_3)e_*^3
-\tilde{T}(J_0e_2)e_*^2-\tilde{T}(J_0e_1)e_*^1= -\tilde{T}J_0.
\end{eqnarray*}
Consequently,
$4\pi (c_1(M,J_1)-c_1(M,J_0))\wedge\omega_0 
= div ((\tilde{T}J_0)^{\sharp})Vol_M
+ d(G \wedge \omega_0).$ 
Integration over $M$ and Stokes on the second term 
leads to the first equality of Theorem 1.2.
Assume $\lnab{}^{E_{J_0}}\tilde{H}$ is anti-$J_0$-complex, that is
$\lnab{J_0X}^{E_{J_0}}\tilde{H}=-J_{0}\lnab{X}^{E_{J_0}}\tilde{H}$.
Then, 
\[  (1+\cos\theta)\tilde{T}(J_0(X))=
\langle\lnab{J_0X}\tilde{H},J_0
\tilde{H}\rangle= \langle\lnab{J_0X}^{E_{J_0}}\tilde{H},J_0
\tilde{H}\rangle=-\langle\lnab{X}\tilde{H},\tilde{H}\rangle.\]
We have $\tilde{H}=\sin\theta H$ and $\langle \lnab{X} H,H\rangle=0$.
Thus, $\langle \lnab{X}\tilde{H},\tilde{H}\rangle=
\ha d\sin^2\theta(X),$
and so,~
$div( (\tilde{T}J_0)^{\sharp})
= -\frac{1}{2}div\LA{(}\frac{\nabla\sin^2\theta}
{(1+\cos\theta)}\LA{)}$. Now $~\frac{\nabla\sin^2\theta}{(1+\cos\theta)}
=\frac{2\sin\theta\cos\theta\, \nabla \theta}{1+\cos\theta}~$, and
\[\begin{array}{l}
\frac{2\sin\theta\cos\theta}{1+\cos\theta}=
\frac{4\sin\frac{\theta}{2}\cos\frac{\theta}{2}\cos\theta}
{2\cos^2\frac{\theta}{2}}=2\tan\frac{\theta}{2}\cos\theta=
(2\sin\frac{\theta}{2}\cos\frac{\theta}{2})\frac{\cos\theta}
{\cos^2\frac{\theta}{2}}=
(2\sin\frac{\theta}{2}\cos\frac{\theta}{2})\LA{(}2-
\frac{1}{\cos^2\frac{\theta}{2}}\LA{)}\\
~~~~~~~~~~~~=\LA{(}4\sin\frac{\theta}{2}\cos\frac{\theta}{2}-2\frac{
\sin\frac{\theta}{2}}{\cos\frac{\theta}{2}}\LA{)}
=\La{(}2\sin\theta -2\tan\frac{\theta}{2}\La{)}
=\LA{(}2\sin\theta-\frac{2\sin\frac{\theta}{2}}{\cos
 \frac{\theta}{2}}\LA{)}
\end{array}\]
That is $\frac{\nabla\sin^2\theta}{(1+\cos\theta)}=\nabla\La{(}\!
-2\cos\theta+ 4\log |\cos\frac{\theta}{2}|\La{)},$ 
where $\theta\in [0,\pi]$, obtaining 
\[
div \La{(}(\tilde{T}J_0)^{\sharp}\La{)}
=-\frac{1}{2}div\LA{(}\frac{\nabla\sin^2\theta}
{(1+\cos\theta)}\LA{)}
=
\Delta \La{(}\cos\theta -2\log (\cos\frac{\theta}{2})\La{)}.
\]
Since $\log(\cos\frac{\theta}{2})=\ha \log\left(\frac{1+\cos\theta}{2}\right)
=\ha(\log(1+\cos\theta)-\log(2))$,
integration over $M$ of $div \La{(}(\tilde{T}J_0)^{\sharp}\La{)}$
 and Stokes on the first term $ \Delta \cos\theta$ gives
 formula (1.4) in the Theorem.\\[3mm]
Now assume $\Sigma$ is a finite disjoint union of connected compact and 
orientable  submanifolds $\Sigma_i$ of dimension $d_i\leq 2$. 
Then for a sufficiently small
tubular neighbourhood of $\Sigma$, $V(\Sigma)=\bigcup_i V(\Sigma_i)$, 
 $\partial (M\sim V(\Sigma))=
\cup_i-\partial V(\Sigma_i)$. Below we  specify such tubular neighbourhood.
For notational simplicity we assume
$\Sigma=\Sigma_i$. 
\\[3mm]
 For each $ r>0$, let $G_r$ and $C_r$ be the subsets  of the total space 
$N\Sigma$ of the normal bundle $T\Sigma^{\bot}$ of $\Sigma$ in $M$
given by
\[G_r=\{(p,w): p\in \Sigma, w\in T_p\Sigma^{\bot}, \|w\|\leq r\},
~~~~~~C_r=\{(p,w): p\in \Sigma, w\in T_p\Sigma^{\bot}, \|w\|=r\}\]
and for $r$ sufficiently small (say $r\leq r_0\leq 1$)
the exponential map of
$M$ restricted to $G_r$, denoted by 
$exp:G_r\ra M$, $exp(p,w)=exp_p(w)$  defines a diffeomorphism onto
$V(\Sigma,r)=\{q\in M: d(q,\Sigma)\leq r\}$ and
$C(\Sigma,r)=\{q\in M: d(q,\Sigma)=r\}$ is its boundary.
Let $S(p,r)$ denote the sphere of radius $r$ in
$T_p\Sigma^{\bot}$,
 $\sigma(q)=d(q,\Sigma)$, and  for each  $w\in T_p\Sigma^{\bot}$,
$\gamma_{(p,w)}(r)=exp_p(rw)$ is 
the  geodesic normal to
$\Sigma$, starting at $p$ with initial velocity $w$.
Thus, $s(p,w):=\sigma(exp(p,w))=\|w\|$, is just the Euclidean norm
in $T_p\Sigma^{\bot}$. Since $N\Sigma$ is the total space of a Riemannian
vector bundle, then it has a natural Riemannian structure
such that $\pi:N\Sigma\ra \Sigma$ is a Riemannian submersion.
The volume element $Vol_{N\Sigma}$ for such metric satisfies
$Vol_{N\Sigma}(p,w)=Vol_{\Sigma}(p)\wedge ds(p,w)$ and
$Vol_{C_r}(p,w)= Vol_{\Sigma}(p)\wedge Vol_{S(p,r)}(w)$, where $r=\|w\|$.
For each $u\in S(p,1)$ define $\vartheta_u(r)= \langle Vol_{N\Sigma}(p,ru),
exp^{*}Vol_M(p,ru)\rangle$. This function $\vartheta_u(r)$ measures the volume
distortion by $exp$ in the direction $u$. It satisfies $\vartheta_u(0)=1$
and a Riccati differential equation.
We recall the following (see \cite{[Gr]}):
(1) $\nu(q)=\nabla \sigma(q)$
is the unit outward of $C(\Sigma,r)$, 
(2) $\nu(\gamma_{(p,u)}(r))=  \gamma'_{(p,u)}(r)$, 
(3) $ds\wedge *ds$ and $d\sigma\wedge *d\sigma$ are
the volume elements of $N\Sigma$ and $M$ respectively.
(4) $*ds$ and $*d\sigma$ are
the volume elements of each hypersurface 
$C_r$ of $N\Sigma$ and $C(\Sigma,r)$ of $M$ respectively, 
and $exp^*(*d\sigma)(p,w)=\vartheta_{\frac{w}{r}}(r)(*ds)(p,w)$, 
where $r=\|w\|$. We have
\begin{eqnarray}
\lefteqn{\int_{M}-\frac{1}{4\pi}
\Delta\log(1+\cos\theta)Vol_M = \lim_{\epsilon\ra 0}
\int_{M\sim V(\Sigma,\epsilon)}-\frac{1}{4\pi}
\Delta\log(1+\cos\theta)Vol_M}\non\\
&=&\lim_{\epsilon\ra 0}\int_{ C(\Sigma,\epsilon)}\frac{1}{4\pi}
\langle \nabla(\log(1+\cos\theta)), \nu\rangle *d\sigma 
=\lim_{\epsilon\ra 0}\int_{ C_\epsilon}\frac{1}{4\pi}
(exp)^*\LA{(}d(\log(1+\cos\theta +1))(\nu) *d\sigma\LA{)} \non\\
&=& \lim_{\epsilon \ra 0}\int_{\Sigma}\frac{1}{4\pi}
\LA{(}\int_{S(p,\epsilon)}d(\log(1+\cos\theta))_{(exp_p(w))}
(\nu)\vartheta_{\frac{w}{\epsilon}}(\epsilon)
d_{S(p,\epsilon)}(w)\LA{)}d_{\Sigma}(p)\non\\
&=& \lim_{\epsilon \ra 0}\int_{\Sigma}\frac{1}{4\pi}
\LA{(}\int_{S(p,1)}\!\!\!\epsilon^{3-d_i}\,\vartheta_{u}(\epsilon)\,
\frac{d}{dr}\La{|}_{r=\epsilon}\log (\phi_{(p,u)}(r))
d_{S(p,1)}(u)\LA{)}d_{\Sigma}(p)
\end{eqnarray}
where $\phi_{(p,u)}(r)=(1+\cos\theta)(\gamma_{(p,u)}(r))$.
The function $\kappa(p,u)=$ the order of the zero $r=0$ of
$\phi_{(p,u)}(r)$, is
defined in the total space of $N^1\Sigma$
 with values in $[1,+\infty]$ and
satisfies $\kappa(p', u')\leq \kappa(p,u)$ for all $(p',u')$ in a 
neighbourhood of $(p,u)$. In particular it is upper semi-continuous, 
and so it is mensurable. 
Now we have for $k=\kappa(p,u)<+\infty$ and $\forall 0<\epsilon<r_0$ 
\begin{eqnarray}
0~ \neq~ \phi_{(p,u)}(\epsilon) &=&\epsilon^{k}\LA{(}A(p,u)+
\epsilon B(p,u,\epsilon)\LA{)}
\end{eqnarray}
\begin{equation}
\mbox{where}~~~~
A(p,u)=\frac{1}{k!}\frac{d^{k}}{dr^k}\phi_{(p,u)}(0)\neq 0,~~~~~~
B(p,u,\epsilon)=\int^1_0 \frac{(1-s)^k}{k!}
\frac{d^{k+1}}{dr^{k+1}}\phi_{(p,u)}(s\epsilon)ds~~~~~~~~~
\end{equation}
Thus, \\[-8mm]
\begin{equation}
\epsilon\frac{d}{d\epsilon}\log (\phi_{(p,u)}(\epsilon))
=\epsilon\frac{\frac{d}{d\epsilon}\phi_{(p,u)}(\epsilon)}
{\phi_{(p,u)}(\epsilon)}
= \kappa(p,u)+ \epsilon \frac{B(p,u,\epsilon)
+\epsilon \frac{d}{d\epsilon}B(p,u,\epsilon)}{A(p,u)+\epsilon
B(p,u,\epsilon)}.
\end{equation}
Assume that for almost all $(u,p)\in N^1\Sigma$, $\kappa(p,u)
<+\infty$. Then
\begin{equation}
(9.10)=\lim_{\epsilon\ra 0}\frac{1}{4\pi}\int_{\Sigma}\int_{S(p,1)}
\!\!\!\!\!\!\!\!\!\!\!\epsilon^{2-d_i}\vartheta_u(\epsilon)\left(\kappa(p,u)+
\frac{\epsilon\La{(}B(p,u,\epsilon)\!+\!\epsilon\frac{d}{d\epsilon}
B(p,u,\epsilon)\La{)}}
{A(p,u)+\epsilon B(p,u,\epsilon)}\right)d_{S(p,1)}d_{\Sigma}(p)
\end{equation}
and for almost all $(u,p)$\\[-5mm]
\[\lim_{\epsilon\ra 0}\vartheta_u(\epsilon)\LA{(}\kappa(p,u)+
\frac{\epsilon(B(p,u,\epsilon)+\epsilon\frac{d}{d\epsilon}B(p,u,\epsilon)
\La{)}}{A(p,u)+\epsilon B(p,u,\epsilon)}\LA{)}=\vartheta_u(0)\kappa(p,u)=
\kappa(p,u).\]
In order to interchange $\lim_{\epsilon \ra 0}$ with $\int_{\Sigma}
\int_{S(p,1)}$ we need  some conditions to be satisfied. 
Under the assumption of controlled zero set, and if $d_i\leq 2$, we may
apply the dominate convergence theorem  to (9.14), and
only the terms with $d_i=2$ do not vanish. The last inequality
follows immediately from (3.7). \qed\\[-3mm]
\begin{Pp} Assume $\Sigma=\cup \Sigma_i$ as in Theorem 1.2. Then
$(1+\cos\theta)$ has controlled zero set if $(1)$,  or $(2)$ holds:\\[2mm]
$(1)$~ $\kappa\in L^1(N^1\Sigma)$ and there exist a constant
$r_0>0$, and non-negative functions $h, d$ on $N^1\Sigma$ 
with $h/d \in L^1(N^1\Sigma)$ 
such that for almost all $(p,u)\!\in\! N^1\Sigma\,, $
$|A(p,u)|\geq  \sup_{r\leq r_0}
|\frac{r}{k!}\frac{d^{k+1}}{d r^{k+1}}\phi_{(p,u)}(r)|+d(p,u) $ 
and $\sup_{r\leq r_0}\LA{\{}|\frac{r}{k!}\frac{d^{k+1}}{d r^{k+1}}
\phi_{(p,u)}(r)|,|\frac{r^2}{k!}\frac{d^{k+2}}{d r^{k+2}}
\phi_{(p,u)}(r)|\LA{\}}\leq  h(p,u)$,
where $k=\kappa(p,u)$.\\[2mm]
$(2)$~~For all $(p,u)$,  $\phi_{(p,u)}(r)$ is polynomial on 
$r\in [0,r_0]$ with coefficient $A(p,u)$ uniformly
bounded away from $0$, i.e. $|A(p,u)|\geq c$ $\forall (p,u)$,
where $c>0$ is constant.
\end{Pp}
\em Proof. \em (1)~For a.e.$(p,u)$, and for  $0<\epsilon < r_0$,
$|A(p,u)+\epsilon B(p,u,\epsilon)|\geq 
|A(p,u)|-\epsilon |B(p,u,\epsilon)| $ $
\geq |A(p,u)|- sup_{r\leq r_0}|
\frac{r}{k!}\frac{d^{k+1}}{d r^{k+1}}\phi_{(p,u)}(r)|\geq d(p,u).$
Moreover $\epsilon |B(p,u,\epsilon)+\epsilon
\frac{d}{d\epsilon}B(p,u,\epsilon)| \leq h(u,p).$
Therefore $|(9.13)|\leq \kappa(p,u) +  \frac{h(p,u)}{d(p,u)}$, 
what implies that 
$(1+\cos\theta)$ has a controlled zero set.
(2) Let $r_0>0$ sufficiently small s.t. all geodesics
$\gamma_{(p,u)}(\epsilon)$ are defined for $0\leq \epsilon <r_0$.
In this case $\kappa$ takes only finite values, and since
it is u.s.c., it has a maximum $k_0$ in $N\Sigma^1$.
Then we may take $C:= \sum_{\mu\leq k_0}max_{(p,u)\in N\Sigma^1,
\epsilon \leq \epsilon_0}
|\frac{1}{(\mu+1)!}\frac{d^{\mu +1}}{d r^{\mu +1}}\phi_{(p,u)}(\epsilon)|$. 
Then we can find $\epsilon_0>0$ s.t. 
$|A(p,u)+\epsilon B(p,u,\epsilon)|\geq
\frac{c}{2}$ for all $\epsilon\leq \epsilon_0$ and all $(p,u)$.
\qed \\[-3mm]
\begin{Cr} If $J_0$ and $J_1$  
 satisfy the conditions of the last part of Theorem 1.2, then:\\[1mm]
$(1)$~ $\La{(}\int_M \|c_1(M,J_1)\|^2Vol_M\La{)}^{\frac{1}{2}}\geq 
\frac{1}{4\sqrt{2}\pi ({Vol(M)})^{\frac{1}{2}}}
\int_{M}(s^M+\ha \|\lnab{}\omega_0\|^2)Vol_M$, with equality iff 
$c_1(M,J_1)$ $
=f\omega_0$ for some smooth function $f:M\ra \R{}$.\\[1mm]
$(2)$~ If, $c_1(M,J_1)$ is $L^2$-orthogonal to $ \omega_0$ 
 then $\int_M s^M Vol_M\leq 0$ and
$\left|\int_M s^M_{J_0}Vol_M\right|
\leq -\int_M s^M Vol_M$. If $\int_M s^M Vol_M=0$ then
$J_0$ is K\"{a}hler.
\end{Cr}
\em Proof. \em  From last inequalities of Theorem 1.2 we clearly obtain (1). 
Moreover,  we have
$0\geq \int_M (s^M+\ha\|\lnab{}\omega_0\|^2)Vol_M$.
In particular $\int_M s^M Vol_M\leq 0$,
and using the equality (3.6) , and   that $s^M\leq s_J^M$,
 we obtain the proof of (2).\qed\\[5mm]
\em Remark $5$. \em  If $M$ is a surface and $f:M\ra \R{}$ is a function
of absolute value type at a zero $p$ (\cite{[E-G-T]}), i.e.\ locally 
$|f(z)|=|z|^k\psi(z)$
where $z$ is an isothermal coordinate around $p$,
and $\psi(z)>0$, then $\kappa(p,u)
= k$ forall $u$.  We also say that $p$ is a zero of homogeneous order 
$k$.\\[4mm]
\em Remark $6$. \em  Assume $J_0$ is almost K\"{a}hler. 
Note that if $J_0$ is K\"{a}hler, 
$E_{J_0}$ is a parallel subbundle of $\bigwedge^2_+TM$.
If $\tilde{H}$ is a non-zero parallel section of $\bigwedge^2_+TM$ then 
by lemma 3.1 $M$ is hyper-K\"{a}hler, with respect to $(J_0, J_1=H)$.
If $\tilde{H}$ is a parallel section of $E_{J_0}$, then
$J_1$ is K\"{a}hler iff $J_0$ is so.\\[5mm] 
\em Proof of Proposition 1.2. \em  
Since $J_0$ is K\"{a}hler, $\omega_0$ is parallel and so harmonic,
and the Weitzenb\"{o}ck operator satisfies $A(\omega_0)=0$ (see (3.2)).
Since $J_1$ is almost K\"{a}hler, then $\omega_1$ is harmonic. 
Polarization
of the Weitzenb\"{o}ck formula (second equality of (3.3))
and using the fact that $A$ is a
symmetric operator, gives \\[-4mm]
\[
\Delta^+(\langle \omega_0,\omega_1\rangle )=
-\langle \Delta \omega_0, \omega_1\rangle - \langle \omega_0,\Delta \omega_1
\rangle +2\langle \lnab{}\omega_0, \lnab{}\omega_1\rangle
+\langle A(\omega_0),\omega_1\rangle +
\langle \omega_0, A(\omega_1)\rangle= 0
\]
that is, $\cos\theta$ is an harmonic function. Since $M$ is closed
$\cos\theta$ is constant. If $\cos\theta\neq \pm 1$
Prop. 9.1 gives the homotopy, and
orthonormalization of Gram-Schmidt of $\omega_0$, $\omega_1$
gives two self-dual 2-forms $\omega_0$  and $\omega_1'$  that correspond
to two $g$-orthogonal complex structures  $J_0$, K\"{a}hler, and
 $J_1$, almost K\"{a}hler, that anti-commute. Lemma 3.1
proves Proposition 1.2. \qed
\section{Applications}
\setcounter{Th}{0}
\setcounter{Pp}{0}
\setcounter{Cr} {0}
\setcounter{Lm} {0}
\setcounter{equation} {0}
\baselineskip .5cm
\subsection{Surfaces}
(1) ~Let $(M, g)$ be a 
compact Riemann surface
of Gauss curvature $K$. If $\xi: U\subset \Co\ra M$ is a conformal coordinate
with $g=\lambda^2<,>$ where $\lambda^2$ is the conformal factor, then
$K= -\ha \Delta \log \lambda^2$ 
and the Euler number of $M$ 
is ${\cal X}(M)=\frac{1}{2\pi}\int_M K Vol_g.$
If we take a conformal change of the metric $\hat{g}=hg$ where $h$ is a 
non-negative smooth real function, then  $Vol_{\hat{g}}=hVol_g$, and
$\xi$ is still a conformal coordinate
for $\hat{M}=(M,\hat{g})$ with conformal factor $h\lambda^2$. Hence
\[\hat{K}=-\ha\hat{\Delta}\log  h\lambda^2=-\ha( h^{-1}\Delta\log h
+h^{-1}\Delta\log \lambda^2)=-\ha h^{-1}\Delta\log h
+h^{-1}K\\[-2mm]\]
\begin{equation}
{\cal X}(\hat{M})={\cal X}(M) 
-\frac{1}{2\pi}\int_M\ha \Delta\log h \,Vol_g
\end{equation}
If on a isothermal coordinate of  each zero or infinity of
$h$, $h$ is of the form $h(z)=|z|^k\psi(z)$ 
where  $\psi\neq 0,\infty$ is a smooth function and $k\in \Z$,  $h$ is 
said to be of \em absolute value type \em at those points (see remark 5). 
In this case
\begin{equation}
-\frac{1}{2\pi}\int_M \Delta\log h Vol_g=\mbox{~number~of~ zeros~of~}h
-\mbox{~number~of~infinities~of~}h
\end{equation} 
where the zeros or infinities of
$h$ are counting including
multiplicity (see [11]).\\[3mm]
(2)~
Consider $\sigma: S^2\ra \R{2}$ the stereographic projection,
$\phi(x_0,x_1,x_2)=\omega=\frac{1}{1-x_0}(x_1,x_2)$, then 
$\sigma$ defines an isometry between
$S^2\sim {(1,0,0)} $ with
the metric $h g_S$ and  $\R{2}$ with the Euclidean metric $g_0$, where
$h=\frac{1}{(1-x_0)^2}$ and $g_S$ is the usual metric of $S^2$. 
Since $(\R{2}, g_0)$ is a 
flat space, its  Euler form is a zero form, and so
$0={\cal X}(\R{2})= {\cal X}(S^2)- \ha(\mbox{~number~of~infinities~of~}h).$ 
Note that $\frac{1}{h}(x_0,x_1,x_3)= (1-x_0)^2$ has a zero of order
4 at $(1,0,0)$ as a  function on $S^2$. 
We could also consider $\sigma$ as an isometry between 
$(S^2, g_S)$ and 
$\R{2}$ with the metric
$\tilde{h}g_0$  where $\tilde{h}=\frac{4}{(1+ \|w\|^2)^2}$
and observe that $\infty$ is a zero of 4-order of $\tilde{h}$. 
Therefore $ {\cal X}(S^2)=2$.\\[3mm]
(3)~ Let $F:M^2\ra N^4$ be a minimal immersed surface into
a K\"{a}hler complex surface $(N,J,g)$ and denote by
$\theta$ the K\"{a}hler angle of $F$, that is $F^*\omega=\cos\theta Vol_M$
where $\omega(X,Y)=g(JX,Y)$. Let $\hat{g}=\sin^2\theta g_M$
where $g_M$ is the induced metric on $M$ by $F$.
The bundle map $\Phi:TM\ra NM$, $\Phi(X)=(JX)^{\bot}$, 
where $NM$ is the normal bundle, is a conformal bundle
map with $\|\Phi(X)\|^2=\sin^2\theta\|X\|^2$.  If $T'$ is the torsion of 
the $\hat{g}$-Riemannian connection $\lnabp{}$ of $M$ that makes $\Phi$
a parallel isometry from $(TM,\hat{g},\lnabp{}\, )$ to $NM$ with its usual
metric and connection $\lnabo{}$, then (see \cite{[S-V]}(3.2))
\[\Phi(T'(X,Y))= \cos\theta (\lnab{}dF(\Jw X,Y)-\lnab{}dF(X,\Jw Y))\]
where $\Jw$ is a $g_M$-orthogonal complex structure on $M$. Minimality
of $F$ is equivalent to the vanishing of  the r.h.s.\  of the above equation.
This means that $T'\!=\!0$, and so $\lnabp{}\,$ is the Levi-Civita connection
for $\hat{M}=(M,\hat{g})$. 
Moreover, one can prove that if $F$ is not a complex
submanifold,  minimality of $F$
implies  $\sin^2\theta$ is of absolute value type at its zeros and so zeros
are isolated and of finite order.
Taking into consideration that $\Phi$ is 
anti-holomorphic for the $g$-orthogonal
complex structures on $TM$ and $NM$ defined by the usual orientation
(and so it is a reverse orientation bundle map), (10.1) and (10.2) gives
\[-{\cal X}(NM)={\cal X}(R')={\cal X}(\hat{M})={\cal X}(M)
+\sum_{p\in {\cal C}^+} order(p)+\sum_{p\in {\cal C}^-} order(p)\]
where ${\cal C}^+$ and ${\cal C}^-$ are respectively the sets of complex and
anti-complex points, zeros of $\sin^2\theta=(1-\cos\theta)(1+\cos\theta)$.
This formula was obtained by \cite{[E-G-T]}, \cite{[We]} and \cite{[Wo]}, 
that we give here a short
proof, using a singular connection.
\subsection{Minimal 4-manifolds}
Let $F:M\ra N$ be a  4-dimensional immersed
Cayley submanifold $M$ of a Calabi-Yau manifold $N$ of complex 
dimension 4. This means that $M$ is calibrated by one of the
Cayley calibrations of $N$ (see e.g.\ \cite{[J]}), and this is equivalent 
to $M$ to be minimal and  with equal K\"{a}hler angles.
In this case the morphism $\Phi:TM \ra NM$,
$\Phi(X)=(JX)^{\bot}$, where $NM$ is the normal bundle, has the
property of being a conformal one, with coefficient of
conformality $\sin^2\theta$, where $\theta$ is the common
K\"{a}hler angle of $M$. For such submanifolds, the complex
points are zeros of finite order of $\sin^2\theta$.
This is treated in \cite{[S-PV]} where using Theorem 1.1
of the present paper one describes the set of complex points of
$M$ as a residue formula involving the Euler (or equivalently,
the Pontrjagin ) numbers of
$M$ and of the normal bundle. The morphism $\Phi$ defines a 
singular metric $\hat{g}=\sin^2\theta g_M$ on $M$, vanishing at
complex points,  and a
connection $\lnabp{}$ that is $\hat{g}$-Riemannian and has
torsion and identifies $p_1(NM)$ with
$p_1(R')$, and  ${\cal X}(NM)$ with ${\cal X}(R')$ as in
the surface case.
\subsection{Isolated pole in dimension 4}
Let $(M,g)$ be a closed Riemannian 4-manifold, $p_0\in M$, and
$h$ a smooth function on $M\sim \{p_0\}$ with $+\infty> h(p)>0$ 
 and $h(p_0)=0$ or $h(p_0)=+\infty$. 
Assume first that  $g$ is flat on
a neighbourhood of  $p_0$, and 
on a  geodesic ball  $B(p_0,\epsilon)\equiv B(0,\epsilon),$
$h(x)=\|x\|^{2k}$ for some 
integer $k$.  Let $\nu$ be the outward
unit of the (geodesic) sphere $S(p_0,\epsilon)$ of $M$ of radius 
$\epsilon$. Let $f=\log h=2k \log \|x\|$.
Then 
$\int_M div_g(P(\nabla f))Vol_M=
\lim_{\epsilon\ra 0}-\int_{S(p_0,\epsilon)}
g\La{(}P(\nabla f)(x), \nu(x)\La{)}
Vol_{S(p_0,\epsilon)}$, where
$P(\nabla f)$ is defined in (1.1).
For each $x\in S(0,\epsilon)$
\begin{equation}
g\La{(}P(\nabla f)(x), \nu(x)\La{)}={(2k)}^2(2k +6)\epsilon^{-3}
\end{equation}
and so for any $\epsilon$ we obtain the  integer
\begin{equation}
-\frac{1}{32\pi^2}\int_{S(p_0,\epsilon)} 
g\LA{(}P(\nabla f),\nu\LA{)}Vol_{S(0,\epsilon)}
=-\frac{1}{2}{k}^2(k +3).
\end{equation}
In the general case, if $(M,g)$ is not flat near $p_0$
 but on a geodesic ball of $p_0$, 
  $f=\log (\tilde{f} h)$ 
where $+\infty>\tilde{f}>0$ is smooth and $h(exp_{p_0}(w))=\|w\|^{2k}$
($e^f$ is of "absolute value type" at $p_0$), using normal coordinates
at $p_0$ we can easily compute $g(P(\nabla f),\nu)$ on a geodesic ball 
of radius  $\epsilon$, and taking the
limit $\epsilon\ra 0$ of the l.h.s.\ of $(10.4)$ we obtain 
  the same  integer on
r.h.s.\ of $(10.4)$.
\section{Quaternionic-Hermitian 8-manifolds}
\setcounter{Th}{0}
\setcounter{Pp}{0}
\setcounter{Cr} {0}
\setcounter{Lm} {0}
\setcounter{equation} {0}
\baselineskip .5cm
If $(M^8,Q,g)$ is an oriented almost quaternionic Hermitian (AQH) 8-manifold
(see [2] for definitions), 
its fundamental form $\Omega$ is a self-dual 4-form
 of norm $\sqrt{\frac{10}{3}}$.
If $M$ is quaternionic K\"{a}hler (QK) then
the cohomology class  $[\Omega]$
of the fundamental 4-form  represents
the first Pontrjagin class of $E$ (see \cite{[K]}, \cite{[Be2]}), 
where $E=span \{J_1,J_2,J_3=J_1J_2\}$ 
is the rank-3 bundle locally defined
by a pair of anti-commuting $g$-orthogonal positive 
almost complex structures
$J_1,J_2$  of $Q$.
If $M$ is only almost QK (AQK) we may compare $[\Omega]$ with $p_1(E)$, and
we may also use the Obata connection to describe characteristic classes,
in a similar way for almost complex structures. 
For each pair  $Q_0$ and $Q_1$  of AQH structures on M, with 
fundamental forms $\Omega_0$ and $\Omega_1$,
we define at each point $p\in M$ an angle $\theta(p)\in [0, \pi]$ as 
$\cos\theta(p) =\frac{3}{10}\langle
\Omega_0(p),\Omega_1(p)\rangle, $ where
 $\langle, \rangle$ is the Hilbert-Schmidt inner product. Then it is defined
a smooth section $\tilde{H}$ of $\bigwedge^4_+(TM)$, orthogonal to
$\Omega_0$ and with $\|\tilde{H}\|^2=\frac{3}{10}\sin^2\theta$, given by
\begin{equation}
\Omega_1=\cos\theta(p) \Omega_0 + \tilde{H}(p).
\end{equation}
If $Q_0$ and $Q_1$
 are both AQK we would like to compare the cohomology classes of $\Omega_i$
 or compare $p_1(E_0)$ with  $p_1(E_1)$, 
in terms of a PDE on $\theta$, and find conditions to conclude:
(a) $\theta$ is constant; or (b) $Q_1$ is QK;  (c) $\Omega_1$ and 
$\Omega_0$ are homotopic; or (d) $E_1$ and $E_2$ are isomorphic.
As in section 9, 
taking $\cos\theta(p,t)$ and
$\tilde{H(p,t)}$ defined for instance as in (9.3) and (9.2), and 
conveniently reparameterizing $t$ near $t=0$, we have a simple answer to (c):
\begin{Pp} If  $Q_0$ and $Q_1$ are two AQH structures with no 
anti-quaternionic points, i.e. points with $\cos\theta(p)=-1$, then a
homotopy $\Omega_t$ of self-dual 4-forms exist from $\Omega_0$ to $\Omega_1$.
Furthermore, if $Q_0$ and $Q_1$ are AQK and  $p\ra\cos\theta (p)$ is constant, 
then the homotopy can be taken by closed self-dual 4-forms.
\end{Pp}
If $\cos\theta(p)$ is constant, then it is also $\cos\theta(p,t)$ defined
by (9.3). Clearly it follows that $\Omega_t=\cos\theta(p,t)\Omega_0
+\tilde{H}(p,t)$ is also closed.
So, existence of anti-quaternionic points can be an obstruction to homotopy, 
but the non-existence of them does not guarantee that we can take 
$\Omega_t$ a fundamental form for some AQH structure, not even locally.
\\[2mm]
Weitzenb\"{o}ck formula (3.2) and second equality of (3.3) also hold for
$4$-forms, where $A=A(R^M):\bigwedge^4TM^*\ra\bigwedge^4TM^*$ is
also symmetric for the Hilbert-Schmidt inner product, as we see next. 
Denote by $\overline{R^M}$ the induced curvature tensor
of $\bigwedge^4TM^*$: if  $\Omega$ is a 4-form, 
$\overline{R^M}(X,Y)\Omega(X_1,\ldots, X_4)=
\sum_{1\leq j\leq 4}-\Omega(X_1,\ldots,R^M(X,Y)X_j,\ldots, X_4).$
Let $(e_{\al})$ be an o.n. basis of $T_pM$.
The Weitzenb\"{o}ck operator applied to $\Omega$
is given by
\[ A(\Omega)(X_1,X_2,X_3,X_4)=\sum_{1\leq \al\leq 8}\sum_{~1\leq k\leq 4}
(-1)^{k} \La{(}\overline{R^M}(e_{\al},X_k)\Omega\La{)} (e_{\al}, X_1,
\ldots, X_{\hat{k}},\ldots, X_4).\]
\begin{Lm} For any distinct $ \al,\be,\ga,\mu,\nu,\sigma,\rho,\eta$:
\[\begin{array}{ll}
\langle A(e^{\al\be\ga\mu}_*), e^{\nu\sigma\rho\eta}_*\rangle &= 
\langle A(e^{\al\be\ga\mu}_*), e^{\al\sigma\rho\eta}_*\rangle = 0\\
\langle A(e^{\al\be\ga\mu}_*), e^{\al\be\rho\eta}_*\rangle &=  \sm{2
R^M(e_\ga,e_\rho,e_\eta,e_\mu)-2R^M(e_{\mu}, e_{\rho},e_{\eta},e_{\ga})
=-2R^M(e_{\ga},e_{\mu},e_{\rho},e_{\eta})}\\
\langle A(e^{\al\be\ga\mu}_*), e^{\al\be\ga\eta}_*\rangle &= \sm{
-2R^M(e_{\al},e_{\mu},e_{\al},e_{\eta})
-2R^M(e_{\be},e_{\mu},e_{\be},e_{\eta})-2R^M(e_{\ga},e_{\mu},e_{\ga},
e_{\eta})
-Ricci^M(e_{\mu},e_{\eta})}
\end{array}\]
and so $\langle A(\Omega),\Omega'\rangle =\langle \Omega, A(\Omega')\rangle$
for any $\Omega,\Omega'$.
\end{Lm}
Note that if $Q_1$  is  AKQ, since $\Omega_1$ is self-dual and closed then 
it is co-closed and so it is an harmonic 4-form. With similar proof
to the one of Proposition 1.2 (section 9), and using Lemma 11.1,
 we have the following conclusion:
\begin{Pp} If $\Omega_0$ is QK, $\Omega_1$ is AQK and $M$ is compact, then
$\theta$ is constant.  If $\cos\theta\neq -1$ then they are
homotopic 4-forms in $H^4_+(M;\R{})$, and  if
$\Omega_1$ is also QK,  $p_1(E_0)=p_1(E_1)$.
\end{Pp}
In this case a homotopy can be given by 
$\Omega_t=\cos(t)\Omega_0+\sin(t)\Omega_2$,
where $\Omega_2=\frac{\Omega_1-\cos\theta\Omega_0}{\sin\theta}$ is
a closed 4-form  orthogonal to $\Omega_0$,
and $\Omega_{t=\theta}=\Omega_1$. If $\Omega_1$ is AQK then so it
is $\Omega_t$.\\[3mm]
Let $(J_1,J_2, J_3)$ be an hyper-K\"{a}hler (HK) structure on $\R{8}$.
This is an o.n. system (of norm $\sqrt{8}$ ) for the Hilbert-Schmidt inner
product. For each $x=(a,b,c)\in \R{3}$ let $J_x=aJ_1+bJ_2+ cJ_3$.
 The elements with $x\in S^2$ are also
called the HK structure of $\R{8}$ determined by
$J_1,J_2$.  Let ${\cal HK}$ be the set of such HK structures 
(for a given orientation of $\R{8}$).
Given an oriented 4-dimensional subspace
$P$ of $\R{8}$ we take $J'_x$, $J_x ''$,  $x\in S^2$,
 the  HK structures on $P$ and on $P^{\bot}$ respectively,
defined by the $g$-orthogonal  positive complex structures. 
Then $J_x(X)=J'_x(X^{\top})+J_x ''(X^{{\bot}})$,
where $X^{\top}$ and $ X^{\bot}$ are respectively the  orthogonal 
projections of $X$ onto $P$ and $P^{\bot}$, defines a HK structure on $\R{8}$.
Note that $J_x$ and its fundamental 4-form $\Omega$ do not depend on the d.o.n. 
basis  we choose for $P$ and for $P^{\bot}$.
Reciprocally,
given $(J_1,J_2,J_3)$ and $X,Y\in \R{8}$ unit vectors, with $H_X\bot H_Y$
then $P=H_X= span \{ X, J_1X,J_2X, J_3X\},\, P^{\bot}=
span \{Y,J_1Y,J_2Y,J_3Y\}\in Gr(4,8)$ and the previous
construction w.r.t $P, P^{\bot}$ recovers $(J_x)_{x\in S^2}$. Let 
\[{\cal T}=\LA{\{} \Omega\in \sm{\bigwedge^4}(\R{8}): \Omega= \frac{1}{6}(
\omega_1\wedge \omega_1 + \omega_2\wedge \omega_2 + \omega_3\wedge
\omega_3): (J_1,J_2,J_3)\in {\cal HK}\LA{\}} \subset 
\sm{\bigwedge^4_+}(\R{8})\]
where $\omega_i(u,v)=g(J_iu,v)$. 
If $M$ is an oriented 8-dimensional manifold
then we define the corresponding  bundles over $M$, ${\cal HK}(M)$
and ${\cal T}(M)$. Now we look for
obstructions for two fundamental 4-forms $\Omega_0$ and $\Omega_1$ on $M$ 
to be homotopic in ${\cal T}(M)$.
Let us assume that on a neighbourhood $U$ of a point of $M$,  
a unit vector field $X$ is defined, and
$E_0= span \{J_1,J_2, J_3\}$ and $E_1= span \{ J'_1, J'_2,J'_3\}$ with
$J_1=J'_1$. On what follows we argue as in \cite{[A-PV-S]}. At each point,
the subspace  $H'_X= span \{X, J'_1 X, J'_2X,J'_3X\}$ is a $J_1$-
complex subspace of real dimension 4. Since $M$ has dimension 8
 there exist a unit vector $Y'\in {H'_X}^{\bot}\cap H_X^{\bot}$. 
One can see from (11.2) below that 
$d=dim(H_X\cap H'_X)=dim(H_X^{\bot}\cap {H'_X}^{\bot})$.
This dimension $d$ is either 2 or 4 (the spaces are $J_1$-complex).
Thus, on $U$, if rank $E_0\cap E_1$ is constantly equal to one, 
then $d$ is constantly equal to 2, and
we may locally smoothly choose $Y'$. Set $x=(1,0,0)$.
Since $TM=H_X\oplus H_{Y'}$, orthogonal sum,
there exist unique smooth maps $y,y',u,u':U\ra \R{3}$ such that
\begin{equation} \begin{array}{l}
B=\{ X, J_{x}X, J_{y}X+J_{y'}Y',J_{x}(J_{y}X+J_{y'}Y')\}\\
B^{\bot}=\{Y',J_{x}Y', J_{u}Y'+J_{u'}X,J_{x}(J_{u}Y'+J_{u'}X)\},
\end{array}
\end{equation}
define a d.o.n. basis of $H'_X$ and $H'_{Y'}$ respectively,
where $\|y\|^2+\|y'\|^2=\|u\|^2+\|u'\|^2=1$, $g(x, y)= g(x, u)=0$, 
$g(x,y')= g(x, u')=0$,  $g(u,y')+ g(y,u')=0$.
Away from $Q_0$-totally-complex points of $H'_X$, 
that is $y(p)\neq 0$, and since $Q_0$-quaternionic 
(and anti-$Q_0$-quaternionic) points
do not exist as well, i.e. $y'(p)\neq 0$, 
we have
\begin{equation}
J_yX+ J_{y'}Y'=J_{\frac{y}{\|y\|}}(\|y\|X+\|y'\|Y)
\mbox{~~~where~~}Y=-J_{\frac{y}{\|y\|}} J_{\frac{y'}{\|y'\|}}Y'\in
H_X^{\bot}.
\end{equation}
Thus, there exist $x_1=x,x_2=\frac{y}{\|y\|},x_3=x_1x_2$ a
d.o.n. basis of $\R{3}$, and $\zeta$ s.t. 
\begin{equation} \begin{array}{l}
B=\{ X, J_{x_1}X, J_{x_2}(cX+sY),J_{x_3}(cX+sY)\}\\
B^{\bot}=\{Y,J_{x_1}Y, J_{x_2}(cY-sX), J_{x_3}(cY-sX)\},
\end{array}
\end{equation}
where $c=\cos\zeta>0$, $s=\sin\zeta>0$,
define a d.o.n. basis of $H'_X$ and $H'_Y$ respectively. Furthermore
$Y\in H_X^{\bot}\cap {H_X'}^{\bot}$ as well. Note that if we make $y\ra 0$,
we are approaching a  $Q_0$-totally-complex point,
we have at $y(p)=0$, $J_yX+ J_{y'}Y'=J_{y'}Y'$, but we may not smoothly  or
continuously extended $x_2$ and $Y$ as
$x_2=y'$ and $Y=Y'$, or smoothly extend $c$ to $0$ and $s$ to $1$ at 
those points.
This ambiguity will imply an ambiguity on the choice of
$e_3,\ldots, e_8$ bellow, and we need to use (11.4) to define the homotopy
in ${\cal T}(M)$.
A similar problem occurs at $Q_0$-quaternionic points of $H'_X$, where 
$H_X\cap H'_X$ becomes 4-dimensional. A smooth
choice of $Y'$ might be impossible.
Set $e_1=X$, $e_2=J_{x_1}X$, $e_3=J_{x_2}X$, $e_4=J_{x_3}X$,
 $e_5=Y$, $e_6=J_{x_1}Y$, $e_7=J_{x_2}Y$, $e_8=J_{x_3}Y$.
Then $\omega_1= "J_{x_1}"=e^{12}+e^{34}+e^{56}+e^{78}$,
$\omega_2= "J_{x_2}" =e^{13}-e^{24}+e^{57}-e^{68}$,
$\omega_3= "J_{x_3}"=e^{14}+e^{23}+e^{58}+e^{67}$. 
We define
$ \omega_4= "J_4"=e^{17}-e^{28}+e^{35}-e^{46}$, 
$~\omega_5= "J_5"=e^{18}+e^{27}+e^{36}+e^{45}$. 
Then we have 
$\omega'_1=\omega_1$, $\omega'_2=c\omega_2+s\omega_4,$ 
$\omega'_3=c\omega_3+s\omega_5$, 
and so locally,
\begin{eqnarray}
E_1 &=& ( J_{x_1}, cJ_{x_2}+sJ_4, cJ_{x_3}+ sJ_5)\\
\Omega_1 &=& c^2 \Omega_0 +\frac{1}{3}sc(\omega_2\wedge \omega_4
+\omega_3\wedge \omega_5) 
+ \frac{s^2}{6}(\omega_1\wedge\omega_1+\omega_4\wedge\omega_4+
\omega_5\wedge \omega_5)\\
\tilde{H} &=& \sm{\frac{s^2}{6}\LA{(}
(\omega_1\wedge \omega_1+\omega_4\wedge \omega_4
+\omega_5\wedge \omega_5)-\frac{7}{15}(\omega_1\wedge \omega_1
+\omega_2\wedge \omega_2+\omega_3\wedge \omega_3)
\LA{)}+ \frac{2}{3}sc(\omega_2\wedge \omega_4)}.~~~~~~~
\end{eqnarray}
We observe that $J_{x_1},
J_{x_2}, J_{x_3}, J_4, J_5$ is an o.n. system in $\bigwedge^2TM$,
and  $J_{x_1}$ anti-commutes with $J_{x_2}$,
$J_{x_3}$, $J_4$, $J_5$, 
and $J_{x_3}=J_{x_1}J_{x_2}$, $J_5= J_{x_1}J_4$, 
$J_{x_2}$ anti-commutes with $J_4$,
 $J_{x_3}$ anti-commutes with $J_5$.
This explains that  (11.5) defines an
HK structure, and it is elementary to verify it spans $E_1$.
We also have $\omega_2\wedge \omega_4=\omega_3\wedge \omega_5$ and 
$\cos\theta= \frac{1}{15}(7+8c^2)\in [\frac{7}{15}, 1]$ is positive,
(and  $c$ can be continuously extended to 
$Q_0$-quaternionic points and to $Q_0$-totally-complex
points). Thus $\cos\theta\neq -1$ everywhere.  
So we may conclude:
\begin{Pp}  If $Q_0$ and $Q_1$ are two almost quaternionic Hermitian 
structures on a 8-manifold $M$ with $rank(E_0\cap E_1)\geq 1$, then $\Omega_1$
and $\Omega_0$ are homotopic in $\bigwedge^4_+(TM)$.
If $E_0\cap E_1$ has constant rank one and $\cos\theta> \frac{7}{15}$, then
$\Omega_0$ and $\Omega_1$ are locally homotopic in ${\cal T}(M)$. 
\end{Pp}
\em Proof. \em Locally we may take for each $t\in [0,1]$
$B_t$ and $B_t^{\bot}$ as in (11.4) replacing $c$ and $s$ by
$ c(t)=\cos(t\zeta), s(t)=\sin(t\zeta)$. Then $\Omega_t
=\frac{1}{6}( \omega_1(t)\wedge \omega_1(t) +\omega_2(t)\wedge \omega_2(t)
+\omega_3(t)\wedge \omega_3(t))\in {\cal T}(M)$
where $J_1(t)=J_{x_2}$,
$J_2(t)=c(t)J_{x_2}+s(t)J_{4}$, $J_3(t)=J_1(t)J_2(t)$. \qed
\begin{Cr} If  $Q_0$ and $Q_1$ are both $QK$,  
then one of the three cases holds: $(a)$~$E_0=E_1$,
$(b)$~ $E_0\cap E_1$ has rank one and  both $Q_0, Q_1$ are locally
HK structures, and $M$ is K\"{a}hler. Furthermore $p_1(E_0)=
p_1(E_1)$.
$(c)$~ $E_0\cap E_1=\{0\}$. 
\end{Cr}
\em Proof. \em  $E_0\cap E_1$ is a
parallel subbundle of $\bigwedge^2(TM)$, and so it has constant rank.
If it has rank $\geq 2$ then it must have rank 3, and so $E_0=E_1$.
So let us assume  it has rank 1. This means that both $E_0$ and $E_1$
have a common compatible parallel complex structure. This implies
$M$ is K\"{a}hler and $E_0$ and $E_1$ are locally HK (see [3]).
\qed\\[2mm]
\em Remark 7. \em If we want to find a homotopy by 
closed $4$-forms and
$Q_0$ or $Q_1$ are not AQK we have to compute
$p_1(E_0)$ and $p_1(E_1)$ directly from local frames $(J_1,J_2,J_3)$
and $(J'_1, J'_2, J'_3)$. Assume we are in the simplest case that $Q_0$ is
HK, $E_0\cap E_1$ has rank $\geq 1$ and $\bigwedge^4_+$ is flat.
We may compare $p_1(E_1)$ with $p_1(E_0)=0$
using a homotopy $p_1(E_t)$, locally defined by $B_t, B_t^{\bot}$ as in
the previous proofs, and using the bundle $\tilde{E}$ over
$M\times [0,1]$, by computing $p_1(\tilde{E})(X, \frac{d}{dt})$ as
we did in section 9, using the remark 4 of section 7. This would lead
to a transgression form with singularities at $Q_0$-totally complex points
and at $Q_0$-quaternionic points, and so to a residue formula.
This is difficult, because $B_t$ are only locally defined, 
and it requires more investigation on
$\bigwedge^4_+$ and ${\cal T}$, namely, under a musical isomorphism,
 by considering the elements 
of ${\cal T}$ as isomorphic and symmetric operators of $\bigwedge^2TM$,
whose eigenvectors are almost complex structures. 
So we will not discuss this problem
in this paper and leave it for later.

\end{document}